\documentclass[12pt,twoside]{amsart}
\usepackage{amssymb}

\nonstopmode

\numberwithin{equation}{section}
\font\tengothic=eufm10 scaled\magstep 1
\font\sevengothic=eufm7 scaled\magstep 1
\newfam\gothicfam
\textfont\gothicfam=\tengothic
\scriptfont\gothicfam=\sevengothic

\DeclareMathOperator{\gin}{gin}
\DeclareMathOperator{\inn}{in}

\DeclareMathOperator{\pnt}{\raise 0.5mm \hbox{\large\bf.}}

\newtheorem{theorem}{Theorem}[section]
\newtheorem{lemma}[theorem]{Lemma}
\newtheorem{proposition}[theorem]{Proposition}
\newtheorem{corollary}[theorem]{Corollary}

\theoremstyle{definition}
\newtheorem{definition}[theorem]{Definition} % \theoremstyle{remark}
\newtheorem{remark}[theorem]{Remark}
\newtheorem{example}[theorem]{Example}

\newtheorem{notation}[theorem]{Notation}

\begin{document}

\title[Componentwise linearity of tetrahedral curves]{On the componentwise
linearity and the minimal free resolution of a tetrahedral curve}
\author{Christopher A.\ Francisco}
\address{Department of Mathematics, University of Missouri, Mathematical
Sciences Building, Columbia, MO 65211, USA}
\email{chrisf@math.missouri.edu}
\author{Juan C.\ Migliore}
\address{Department of Mathematics, University of Notre Dame,
Notre Dame, IN 46556, USA}
\email{Juan.C.Migliore.1@nd.edu}
\author{Uwe Nagel}
\address
{Department of Mathematics, University of Kentucky,
715 Patterson Office Tower, Lexington, KY 40506-0027, USA}
\email{uwenagel@ms.uky.edu}

%\date{\today}
\thanks{Part of the work for this  paper was done while the second author
was sponsored by the National Security Agency  under Grant Number
MDA904-03-1-0071.}

\begin{abstract}
A tetrahedral curve is an unmixed, usually non-reduced, one-dimensional
subscheme of projective 3-space whose homogeneous ideal is the intersection of
powers of the ideals of the six coordinate lines.  The second and third authors
have shown that these curves have very nice combinatorial properties, and they
have made a careful study of the even liaison classes of these 
curves.  We build
on this work by showing that they are ``almost always" componentwise linear,
i.e.\ their homogeneous ideals have the property that for any  $d$, the degree
$d$ component of the ideal generates a new ideal whose minimal free resolution
is linear.  The one type of exception is clearly spelled out and studied as
well.  The main technique is  a careful study of the way that basic
double linkage behaves on tetrahedral curves, and the connection to the
tetrahedral curves that are {\em minimal} in their even liaison classes.  With
this preparation, we also describe the minimal free resolution of a tetrahedral
curve, and in particular we show that in any fixed even liaison class there are
only finitely many tetrahedral curves with linear resolution.  Finally, we
begin the study of the generic initial ideal ($\gin$) of a tetrahedral curve.
We produce the
$\gin$ for arithmetically Cohen-Macaulay tetrahedral curves and for minimal
arithmetically Buchsbaum tetrahedral curves, and we show how to obtain it for
any non-minimal tetrahedral curve in terms of the $\gin$ of the minimal curve
in that even liaison class.

\end{abstract}

%%%%%%%%%%%%%%%%%%%%%%%%%%%%%%%%%%%%%
\maketitle
\tableofcontents
%%%%%%%%%%%%%%%%%%%%%%%%%%%%%%%%%%%%%%%%%%%%

\section{Introduction}

A \emph{tetrahedral curve} is a curve in $\mathbb{P}^3$ defined by an
ideal
$$
I=(a,b)^{a_1} \cap (a,c)^{a_2} \cap (a,d)^{a_3} \cap (b,c)^{a_4} \cap
(b,d)^{a_5} \cap (c,d)^{a_6} \subset k[a,b,c,d].
$$
These ideals are unmixed of codimension two, and their name comes from
the fact that one can view the six lines defined by the ideals of two of
the variables as forming the edges of a tetrahedron. In his unpublished
Ph.D. thesis \cite{schwartau}, Phil Schwartau studied the case in which
$a_2=a_5=0$, giving a characterization of when the curves are
Cohen-Macaulay in terms of the $a_i$ and describing their minimal free
resolutions.   Note that when $a_2=a_5=0$, the remaining four lines of support
form a complete intersection of type (2,2).

The general case of a tetrahedral curve, when $a_2$ and $a_5$ are not
necessarily zero, is studied in \cite{MN7}. There is a straightforward
reduction procedure for tetrahedral curves using basic double linkage.
Starting with a tetrahedral curve, one does a sequence of basic double
links, getting progressively smaller tetrahedral curves and ending with
one of two outcomes. The reduction process could stop with the empty set,
which we will call the \emph{trivial curve}, defined by the 6-tuple
$(0,0,0,0,0,0)$. Alternatively, one might reach a \emph{minimal curve}
that cannot be reduced further. An easy numerical test allows one to
determine when one has reached a minimal curve, leading to a simple
algorithm for the reduction process. Moreover, all the curves in a
reduction sequence are in the same even liaison class.

The resolutions of the minimal tetrahedral curves have a particularly nice
form. The authors of \cite{MN7} find their graded Betti numbers
explicitly and show that the resolutions are all linear. Additionally,
the length of the resolution guarantees that the trivial curve is the
only minimal arithmetically Cohen-Macaulay curve. A consequence of the
form of the minimal free resolution of minimal tetrahedral curves is that
a tetrahedral curve is minimal (in the reduction process) if and only if
it is minimal in its even liaison class. As applications, the authors of
\cite{MN7} give a new proof of Schwartau's result characterizing the
Cohen-Macaulay curves with $a_2=a_5=0$, classify the 6-tuples of minimal,
arithmetically Buchsbaum curves, and explore unobstructedness and the
Hilbert scheme of some tetrahedral curves.

In this paper, much of our work is devoted to determining when the ideal
of a tetrahedral curve is componentwise linear and the consequences of
this characterization. We recall the definition of a componentwise linear
ideal.

\begin{definition}
Let $I$ be a homogeneous ideal, and write $(I_d)$ for the ideal generated
by the degree $d$ elements of $I$. We say that $I$ is \emph{componentwise
linear} if $(I_d)$ has a linear resolution for all $d$.
\end{definition}

Of course, any ideal with a linear resolution is also componentwise
linear. Some other common examples of componentwise linear ideals include
strongly stable ideals, squarefree strongly stable ideals, and the
${\mathbf a}$-stable ideals of \cite{GHP}.

Componentwise linear ideals were introduced in a paper of Herzog and Hibi
\cite{HH}. Initially, a primary motivation for studying componentwise
linear ideals came from combinatorics and the desire to generalize the
notion of having a linear resolution. Eagon and Reiner proved that if
$\Delta$ is a simplicial complex, and $I_\Delta \subset
R=k[x_1,\dots,x_n]$ is its Stanley-Reisner ideal, then $I_\Delta$ has a
linear resolution if and only if the Alexander dual $\Delta^*$ is
Cohen-Macaulay over $k$ \cite{ER}. Componentwise linear ideals help
extend this statement; $I_\Delta$ is componentwise linear if and only if
$\Delta^*$ is sequentially Cohen-Macaulay, a property that requires a
nice filtration on $R/I_{\Delta^*}$ in which the quotients are
Cohen-Macaulay \cite{HH,HRW}.

In addition, componentwise linear ideals have a number of algebraic
properties that make them interesting to study. Herzog and Hibi proved
convenient formulas for the graded Betti numbers of a componentwise
linear ideal $I$ in terms of the Betti numbers of the $(I_d)$ and
$\mathfrak m (I_d)$, where $\mathfrak m$ is the maximal homogeneous
ideal. Moreover, Aramova, Herzog, and Hibi proved that if the
characteristic of $k$ is zero, and $\gin(J)$ is the reverse-lex generic
initial ideal of $J$, then $J$ and $\gin(J)$ have the same graded Betti
numbers if and only if $J$ is componentwise linear. Thus componentwise
linear ideals have the same graded Betti numbers as strongly stable
ideals, so there is a lot of structure in their resolutions.

The origin of this work is a confluence of ideas from two places. A
remark in \cite{MN7} notes that there are a number of linear strands in
the minimal free resolution of the ideal of a tetrahedral curve. We
wanted to find a clear explanation for how these linear strands arise.
Additionally, the main result of \cite{Fran} is that ideals of at most
$n+1$ general fat points in $\mathbb{P}^n$ are componentwise linear. One
can take these ideals to be the intersection of powers of ideals
generated by $n$ of the $n+1$ variables; that is, in $\mathbb{P}^3$, they
have the form
$$
(b,c,d)^{b_1} \cap (a,c,d)^{b_2} \cap (a,b,d)^{b_3} \cap (a,b,c)^{b_4}.
$$

These ideals are similar enough to the ideals of tetrahedral curves that
we wondered if one might be able to prove that some large class of
tetrahedral curves is componentwise linear, and tests using Macaulay 2
\cite{M2} and the MAPLE code from \cite{MN7} suggested many of our
results in the following sections.

The main tool throughout our paper is the reduction process for
tetrahedral curves from \cite{MN7}. We begin our investigation in section
\ref{bdl}  by determining
in Proposition \ref{theorem A} how componentwise linearity persists in a
basic double link. This analysis forms the basis for Theorem \ref{theorem
B} and Corollary \ref{cwl conclusion}, which characterize which ideals of
tetrahedral curves are componentwise linear in terms of the curves to
which they reduce. In the case of Schwartau curves, when $a_2=a_5=0$, we
can say more, proving in Corollary \ref{schwartau curves} that ideals of
Schwartau curves fail to be componentwise linear if and only if
$a_1+a_6=a_3+a_4$, and $a_1$, $a_3$, $a_4$, and $a_6$ are all positive.

As applications of our results on componentwise linearity, we prove a
number of statements about the minimal free resolutions of ideals of
tetrahedral curves. The ideals $J$ that are not componentwise linear are
actually not far from being componentwise linear, which we measure in
Proposition \ref{ginres} by comparing the graded Betti numbers of the
reverse-lex generic initial ideal $\gin(J)$ to those of $J$. One
consequence is Theorem \ref{reg-weight}, which gives an explicit
expression of the regularity of any tetrahedral curve in terms of the
$a_i$. Additionally, Corollary \ref{cor-mfr} describes an easy iterative
procedure for calculating the graded Betti numbers of any tetrahedral curve
from just the
$a_i$ and a knowledge of the graded Betti numbers of the minimal curves from
\cite{MN7}.

In section  \ref{linear
resolutions}, we investigate which tetrahedral curves, in addition to the
minimal ones, have linear resolutions. We characterize the arithmetically
Cohen-Macaulay curves with linear resolutions in Proposition \ref{acm lin
res} and find all the tetrahedral curves with linear resolutions that are
in the even liaison class of two skew lines in Proposition \ref{class of 2
lines}. In addition, we show in Theorem
\ref{finitely many curves} that there are only finitely many tetrahedral curves
with a linear resolution in the even liaison class of a tetrahedral curve that
is not arithmetically Cohen-Macaulay.

Finally, we conclude with some observations about the reverse-lex generic
initial ideal of a tetrahedral curve. The gin is easy to describe in the
arithmetically Cohen-Macaulay case, and we discuss how the gin changes
with a basic double link in the non-arithmetically Cohen-Macaulay case.  In
particular, if we know the gin for a minimal non-arithmetically Cohen-Macaulay
tetrahedral curve then we know it for any tetrahedral curve in the
corresponding even liaison class.  We carry out this program for the
arithmetically Buchsbaum tetrahedral curves.

Throughout, we will often abuse notation and refer to the ideal $I=(a_1,
\dots,a_6)$ or the curve $C=(a_1,\dots,a_6)$ interchangeably.

%%%%%%%%%%%%%%%%%%%%%%%%%%%%%%%%%%%%%%%%%%%%%%%%%%%%%%%%%%%%%%%%%%%%%%

\section{Preliminaries}

 We will denote by $R$
the polynomial ring $k[x_0,x_1,\dots,x_n]$, where $k$ is any field.  We also
denote by ${\mathfrak m}$ the irrelevant ideal $(x_0,x_1,\dots,x_n)$.  Starting
with section \ref{tetra CWL}, though, we will follow \cite{schwartau} and
\cite{MN7} and let $R = k[a,b,c,d]$.

\begin{remark}
When we refer to ``the smaller curve" in the proofs in this paper, we mean
the smaller of the schemes defined by the corresponding ideals, not the
smaller of the ideals.
\end{remark}

\begin{notation}
For a homogeneous ideal $I \subset R = k[x_1,\dots,x_n]$, we let
$I_{\ge d}$ be the ideal generated by all elements of $I$ of degree at least
$d$. Furthermore, $I_d$ will denote the degree $d$ part of $I$, and 
$(I_d)$ will
denote the ideal generated by the degree $d$ part of $I$.
\end{notation}

We begin with a lemma describing how the graded Betti numbers of $I$ and
$I_{\ge d}$ differ.

\begin{lemma}\label{truncatebetti} Let $I$ be a homogeneous ideal in
$S=k[x_1,\dots,x_n]$, and let $d$ be a positive integer. Then for each integer
$r
\ge 0$ and all $i$,
$$
  \beta_{i,i+d+1+r}(I_{\ge d}) = \beta_{i,i+d+1+r}(I).
$$
\end{lemma}

\begin{proof} We have the short exact sequence
$$
0 \longrightarrow I_{\ge d} \longrightarrow I \longrightarrow I/I_{\ge d} \longrightarrow
0.
$$
This induces a long exact sequence in Tor: For all $r \ge 0$,
$$
  \cdots \rightarrow \text{Tor}_{i+1}(I/I_{\ge d},k)_{i+d+1+r} \rightarrow
\text{Tor}_i(I_{\ge d},k)_{i+d+1+r} \rightarrow \text{Tor}_i(I,k)_{i+d+1+r}
\rightarrow
$$
$$
  \text{Tor}_i(I/I_{\ge d},k)_{i+d+1+r} \rightarrow \cdots
$$
is an exact sequence of $k$-vector spaces. Moreover, $I/I_{\ge d}$ has finite
length; it is zero in degree $d$ and higher and has highest degree socle
generator in degree $d-1$. Therefore $I/I_{\ge d}$ has regularity $d-1$,
meaning
$\beta_{i,i+d+r}(I/I_{\ge d})=0$ for all $i$ and all $r \ge 0$.

For $r \ge 0$, because $d+r>d-1$,
$$
  \dim_k \text{Tor}_{i+1}(I/I_{\ge d},k)_{i+d+1+r} = 
\beta_{i+1,i+d+1+r}(I/I_{\ge
d}) = 0.
$$
Similarly,
$$
  \dim_k \text{Tor}_{i}(I/I_{\ge d},k)_{i+d+1+r} = 
\beta_{i,i+d+1+r}(I/I_{\ge d})
= 0.
$$
Consequently, as $k$-vector spaces,
$$
\text{Tor}_i(I_{\ge d},k)_{i+d+1+r} \cong \text{Tor}_i(I,k)_{i+d+1+r}.
$$ Hence their dimensions over $k$ are equal, and thus for all $r \ge 0$,
$$
\beta_{i,i+d+1+r}(I_{\ge d}) = \beta_{i,i+d+1+r}(I).
$$
\end{proof}

A very basic tool used in \cite{MN7} and in this paper is that of {\em basic
double linkage}.  This very simple but powerful construction was introduced by
Lazarsfeld and Rao \cite{LR} to describe the even liaison class of a general
curve in ${\mathbb P}^3$, but it has seen a wealth of generalizations and
applications since then, far too many to list here.  We refer the reader to
\cite{migbook} for some of these, although many more have emerged since
\cite{migbook} was published.  We recall here the codimension two construction
and important facts of basic double linkage, and even this will be a special
case (using a linear form instead of a form of any degree) for the purposes
needed below.  Again, we cite \cite{migbook} for the proofs.  For convenience,
in the result below we denote by $\deg(I)$ the degree of the scheme defined by
$I$.

\begin{theorem} \label{bdl results}
Let $I \subset R$ be a homogeneous ideal, and let $F \in I$ be a homogeneous
polynomial of degree $d$.  Let $L \in R_1$ be a linear form such that 
$L$ is not
a factor of $F$, i.e.\ such that $(L,F)$ is a regular sequence.  Let $J$ be the
ideal $L \cdot I + (F)$.  $J$ is called a {\em basic double link} 
of $I$. Then

\begin{itemize}

\item[(a)] We have an exact sequence
\[
0 \rightarrow R(-d-1) \rightarrow I(-1) \oplus R(-d) \rightarrow J \rightarrow
0
\]
where the first map is given by $C \mapsto (FC, LC)$ and the second is given by
$(A,B) \mapsto LA - FB$.

\item[(b)] $J$ is saturated if and only if $I$ is saturated.

\item[(c)] $J$ is unmixed if and only if $I$ is the saturated ideal of a
codimension two subscheme of ${\mathbb P}^n$.  In this case we have $\deg(J) =
\deg(I) + \deg(F)$.

We assume from now on that $J$ is unmixed.

\item[(d)] $J$ is linked in two steps to $I$.  Hence basic double
linkage preserves the even liaison class of $I$.  In particular, $J$ is
arithmetically Cohen-Macaulay if and only if $I$ is arithmetically
Cohen-Macaulay.  Also, $J$ is locally Cohen-Macaulay and equidimensional if and
only if $I$ is locally Cohen-Macaulay and equidimensional.

\end{itemize}

\end{theorem}

%%%%%%%%%%%%%%%%%%%%%%%%%%%%%%%%%%%%%%%%%%%%%%%%%%%%%%%%%%%%%%%%%%%%%

\section{Basic Double Linkage and Componentwise Linear Ideals} \label{bdl}

In this section we find initial connections between the construction of basic
double linkage in codimension two and componentwise linear ideals.  These will
be important in the subsequent sections.

\begin{proposition} \label{theorem A}
Assume that  $I$ is
componentwise linear.  Let $F \in I$ and let $L$ be a linear form such
that $(L,F)$ is a regular sequence.  Let $J = L \cdot I + (F)$.  Then $J$
is componentwise linear if and only if $F$ is {\em not} a minimal
generator of $I$.
\end{proposition}

\begin{proof}
First assume that $F$ is not a minimal generator of $I$.   Let
$\deg F = e$.  For any $d$, we have that  $J_d = L \cdot I_{d-1} + F \cdot
{\mathfrak m}^{d-e}$, where we make the convention that ${\mathfrak
m}^{d-e} = 0$ for $d < e$.  If $d < e$ then $J_d = L \cdot I_{d-1}$, so
$(J_d)$ has a linear resolution since $(I_{d-1})$ does.

Next we suppose that $d=e$.  We have $J_e$ is spanned by $L \cdot
I_{e-1}$ and $F$.  It follows that $(J_e) = L \cdot (I_{e-1}) + (F)$.
Since $F$ is not a minimal generator of $I$, $F \in (I_{e-1})$.  Hence
the ideal $(J_e)$ arises as a basic double link from $(I_{e-1})$ using
$L$ and $F$.  We  then have from Theorem \ref{bdl results} the exact sequence
\begin{equation} \label{liai seq}
0 \rightarrow R(-e-1) \rightarrow (I_{e-1})(-1) \oplus R(-e) \rightarrow
(J_e) \rightarrow 0.
\end{equation}
The mapping cone then gives a linear resolution for $(J_e)$.

Finally, let $d>e$.  We know from Theorem \ref{bdl results} that we have an
exact sequence
\begin{equation} \label{BDL sequence}
0 \rightarrow R(-e-1) \rightarrow I(-1) \oplus R(-e) \rightarrow J
\rightarrow 0.
\end{equation}
We claim that we now have a short exact sequence
\begin{equation} \label{chris seq}
0 \to {\mathfrak m}^{d-e-1}(-e-1) \to (I_{d-1})(-1) \oplus {\mathfrak
  m}^{d-e}(-e) \to (J_d)
\to 0.
\end{equation}
Indeed, the first map is given by $C \mapsto (FC, LC)$ and the second map
is given by $(A,B) \mapsto LA - FB$.  Because $F \in I_{e}$, the kernel
of the second map is immediately seen to be isomorphic to ${\mathfrak
m}^{d-e-1}$, since $(L,F)$ is a regular sequence and so $L$ and $F$ have
no common factor.  This gives a  diagram
\begin{equation} \label{chris diag}
\begin{array}{cccccccccccccccccc}
&& \vdots &&& \vdots \\
&& \downarrow &&& \downarrow \\
&& R(-d-1)^\bullet && R(-d-1)^\bullet & \oplus & R(-d-1)^\bullet \\
&& \downarrow &&& \downarrow \\
&& R(-d)^\bullet && R(-d)^\bullet & \oplus & R(-d)^\bullet \\
&& \downarrow &&& \downarrow \\
0 & \rightarrow & \mathfrak m^{d-e-1}(-e-1) & \rightarrow & (I_{d-1})(-1)
& \oplus & \mathfrak m^{d-e}(-e) & \rightarrow & (J_d) & \rightarrow & 0
\\
&& \downarrow &&& \downarrow \\
&& 0 &&& 0
\end{array}
\end{equation}
Since $F {\mathfrak m}^{d-e-1} \subset (I_{d-1})$ we get for all $i$
that every minimal $i$-th  syzygy of ${\mathfrak m}^{d-e-1}$ is an
$i$-th syzygy of $(I_{d-1})$, thus a minimal syzygy for degree
reasons.   Hence, the terms in the mapping cone coming from
the leftmost column all get split off, leaving a linear resolution for
$(J_d)$.  This completes one direction of the proof.

Conversely, we assume that $F$ {\em is} a minimal generator of $I$, and
we show that then $J$ is not componentwise linear.  Again suppose $\deg F
= e$.  We have $J = L \cdot I + (F)$.  We will show that $(J_e)$ does not
have a linear resolution.  Note that $J_e $ is again spanned by $L \cdot
I_{e-1}$ and $F$.  Consider the exact sequence
\[
0 \rightarrow K \rightarrow (I_{e-1})(-1) \oplus R(-e) \rightarrow (J_e)
\rightarrow 0
\]
where the second map is given by $(A,B) \mapsto LA - FB$, and $K$ is just
the kernel.  An element of the kernel of this map corresponds to a pair
$(A,B)$ for which $LA = FB$.  An element of $K$ of degree $e$ corresponds
to a pair $(A,\lambda)$, where $A \in I_{e-1}$, $\lambda \in k$, and $LA
= \lambda F$.  But $F \notin (I_{e-1})$ since $F$ is a minimal generator,
so this is impossible.  Hence $K_e = 0$.  An element of $K$ of degree
$e+1$ corresponds to a pair $(A,B)$ where $A \in (I_{e-1})_e$, $B$ is a
linear form, and $LA = FB$.  But $L$ and $F$ have no common factor, so up
to scalar multiple we have $L = B$ and $A = F$.  But $F$ is a minimal
generator of $I$, so again $F \notin (I_{e-1})$.  We thus have that also
$K_{e+1} = 0$.  But this means that $K$ has generators in degree $\geq
e+2$.  Since by hypothesis $(I_{e-1})$ has a linear resolution, the
mapping cone gives a resolution for $(J_e)$ that cannot be linear.
\end{proof}

\begin{corollary} \label{where resol is not lin}
Assume that $I$ is componentwise linear, and let $J = L\cdot I + (F)$,
where $F \in I$ and $(F,L)$ is a regular sequence.  Assume further that
$F$ is a minimal generator of $I$ of degree $e$.  Then $(J_d)$ has a
linear resolution if and only if $d \neq e$.
\end{corollary}

\begin{proof}
We know from Theorem \ref{theorem A} that there is at least one $d$ for
which $(J_d)$ does not have a linear resolution.  We have  seen, in fact,
that
$(J_e)$ does not have a linear resolution, proving one direction here.
So we have only to show that if $d \neq e$ then $(J_d)$ has a linear
resolution.  The proof is very similar to that of Theorem \ref{theorem
A}.  If $d < e$ then $J_d = L \cdot I_{d-1}$, and the linearity is
clear.  If $d > e$ then (\ref{chris seq}) and (\ref{chris diag}) continue
to hold, and the linearity of the resolution is proved in the same way.
\end{proof}

\begin{corollary} \label{first cor}
Let $F \in I$ and let
$L$ be a linear form such that $(L,F)$ is a regular sequence.  Let $J = L
\cdot I + (F)$.  Assume that $F$ is not a minimal generator of $I$.  Then

\begin{itemize}
\item[(a)] $J$ is componentwise linear if and only if $I$ is componentwise
linear.

\item[(b)] If $I$ has a linear resolution and generators of degree $e-1$
then $J$ has a linear resolution if and only if $\deg F = e$.

\item[(c)] If $J$ has a linear resolution then so does $I$.
\end{itemize}
\end{corollary}

\begin{proof}
Part (a) is a subtle variation of Theorem \ref{theorem A} which will,
nevertheless, prove useful.  If $I$ is componentwise linear then we have
already proved the result in Theorem \ref{theorem A}.  So we must assume
that $J$ is componentwise linear.  The proof is almost identical to that
given in Theorem \ref{theorem A}.  We still get the diagram (\ref{chris
diag}) if $d > e$, and so we only have to observe that in order for the
resulting
resolution for $J_d$ to be linear, we must have the resolution for
$I_{d-1}$ be linear as well.

If $d = e$, then we argue similarly by using the sequence
(\ref{liai seq}).

Parts (b) and (c) follow using similar arguments, using the sequence
(\ref{BDL sequence}).
\end{proof}

%%%%%%%%%%%%%%%%%%%%%%%%%%%%%%%%%%%%%%%%%%%%%%%%%%%%%%%%%%%%%%%%%%%%%%

\section{When is a Tetrahedral Curve Componentwise Linear?} \label{tetra CWL}

We now apply the results of the preceding section to tetrahedral curves.
 From now on $R$ will denote the ring $k[a,b,c,d]$.  We first recall some
basic results from \cite{MN7}.

\begin{proposition}[\cite{MN7} Proposition 3.1] \label{MN7 basic result}
Let $J = (a,b)^{a_1} \cap (a,c)^{a_2} \cap (a,d)^{a_3} \cap
(b,c)^{a_4} \cap (b,d)^{a_5} \cap (c,d)^{a_6}$ where not all exponents $a_i$
are zero.  Consider the following systems of inequalities:
\[
\begin{array}{ccc}
\begin{array}{rrcl}
(A): & a_1+a_2 & \geq & a_4,\\
& a_1+a_3 & \geq & a_5, \\
& a_2+a_3 & \geq & a_6
\end{array}
&&
\begin{array}{rrcl}
(B): & a_1+a_4 & \geq & a_2,\\
& a_1+a_5 & \geq & a_3, \\
& a_4+a_5 & \geq & a_6
\end{array}
\\
\begin{array}{rrcl}
(C): & a_2+a_4 & \geq & a_1,\\
& a_2+a_6 & \geq & a_3, \\
& a_4+a_6 & \geq & a_5
\end{array}
& \hbox{\hskip 1cm} &
\begin{array}{rrcl}
(D): & a_3+a_5 & \geq & a_1,\\
& a_3+a_6 & \geq & a_2, \\
& a_5+a_6 & \geq & a_4.
\end{array}
\end{array}
\]
For $1 \leq i \leq 6$ let $a_i' = \max \{ 0, a_i -1 \}$.  Then we have
\begin{itemize}
\item[(i)] $(A) \Leftrightarrow$ $J$ is a basic double link of
\[
(a,b)^{a_1'} \cap (a,c)^{a_2'} \cap
(a,d)^{a_3'} \cap (b,c)^{a_4} \cap (b,d)^{a_5} \cap (c,d)^{a_6}
\]
using $F = b^{a_1}c^{a_2}d^{a_3}$ and $G = a$.
\item[(ii)] $(B) \Leftrightarrow$ $J$ is a basic double link of
\[
(a,b)^{a_1'} \cap (a,c)^{a_2} \cap (a,d)^{a_3} \cap (b,c)^{a_4'} \cap
(b,d)^{a_5'} \cap (c,d)^{a_6}.
\]
using $F = a^{a_1}c^{a_4}d^{a_5}$ and $G = b$.
\item[(iii)] $(C) \Leftrightarrow$ $J$ is a basic double link of
\[
(a,b)^{a_1} \cap (a,c)^{a_2'} \cap (a,d)^{a_3} \cap (b,c)^{a_4'} \cap
(b,d)^{a_5} \cap (c,d)^{a_6'}.
\]
using $F = a^{a_2}b^{a_4}d^{a_6}$ and $G = c$.
\item[(iv)] $(D) \Leftrightarrow$ $J$ is a basic double link of
\[
(a,b)^{a_1} \cap (a,c)^{a_2} \cap (a,d)^{a_3'} \cap (b,c)^{a_4} \cap
(b,d)^{a_5'} \cap (c,d)^{a_6'}.
\]
using $F = a^{a_3}b^{a_5}c^{a_6}$ and $G = d$.
\end{itemize}

\end{proposition}

\begin{remark} \label{remark-numerical}
There is no known numerical criterion that characterizes whether a tetrahedral curve is arithmetically Cohen-Macaulay in terms of the $a_i$ (except for the unpublished result of Schwartau in the case where $a_2 = a_5 = 0 $ (cf. \cite{MN7}, Theorem 5.3)). However, a necessary and sufficient condition for a tetrahedral curve to be arithmetically Cohen-Macaulay is for there to exist a sequence of reductions of the form given in Proposition~\ref{MN7 basic result} down to a complete intersection (and ultimately to the trivial curve) \cite{MN7}. We can carry out this reduction process by sequentially reducing facets of maximal weight; see Lemma~\ref{reduce max wt lemma} and Example~\ref{reduce-example}.
\end{remark}

\begin{definition}[\cite{MN7} Theorem 5.1] \label{def of min}
A non arithmetically Cohen-Macaulay tetrahedral curve $C$ is {\em minimal} if
either of the following equivalent conditions holds:

\begin{itemize}
\item[(a)] The ideal $I_C$ does not admit any reduction of the type given in
parts (A) to (D) of Proposition \ref{MN7 basic result};

\item[(b)] $C$ is minimal in its even liaison class (cf.\ \cite{migbook}).
\end{itemize}
\end{definition}

\begin{corollary}[\cite{MN7} Corollary 3.5] \label{s-min}
Consider a tetrahedral curve $C = (a_1,a_2,a_3,a_4,a_5,a_6)$ where
not all $a_i$
are 0.  Assume without loss of generality that $a_6 = \max \{ a_1,\dots,a_6
\}$.    Then $C$ is minimal if and only if
\[
\begin{array}{rcl}
a_1 & > & \max \{ a_3 + a_5,a_2+a_4\} \hbox{ and} \\
a_6 & > & \max \{ a_4+a_5,a_2+a_3 \}
\end{array}
\]
\end{corollary}

\begin{theorem}[\cite{MN7} Theorem 4.2] \label{thm-res}
Every non-trivial {\em minimal} tetrahedral curve has a linear minimal
free resolution.

More precisely, if the curve $C$ is defined by $(a_1,a_2,a_3,a_4,a_5,a_6)$
and
$a_6 = \max \{a_i\} > 0$ then its minimal free resolution has the form
\begin{equation*}
     0 \to R^{\beta_3}(-a_1-a_6-2)  \to R^{\beta_2}(-a_1-a_6-1)  \to
R^{\beta_1}(-a_1-a_6) \to I_C \to 0
\end{equation*}
where
\begin{eqnarray*}
\beta_1 & = & (a_1 + 1) (a_6 + 1) - \sum_{i=2}^5 \frac{a_i (a_i + 1)}{2} \\
\beta_2 & = & 2 a_1 a_6 + a_1 + a_6  - \sum_{i=2}^5 a_i (a_i + 1) \\
\beta_3 & = & a_1 a_6  - \sum_{i=2}^5 \frac{a_i (a_i + 1)}{2}. \\
\end{eqnarray*}
\end{theorem}

In order to have an (almost) canonical way to reduce to a  minimal
tetrahedral curve, we use facets of maximal weight.  Recall that the {\em
weight of a facet} is the sum of the weights of the edges forming its
boundary.

\begin{lemma}[\cite{MN7} Lemma 3.8] \label{reduce max wt lemma}
Let $C =  (a_1,a_2,a_3,a_4,a_5,a_6)$ be a non-trivial tetrahedral curve.  If
$C$ is not minimal then one can reduce any of its facets of maximal weight.
\end{lemma}

\begin{example} \label{reduce-example}
Consider the curve $(3,3,3,1,2,4)$.  The facets have the following weights:
$a_1+a_2+a_3 = 9$, $a_1+a_4+a_5 = 6$, $a_2+a_4+a_6 = 8$ and $a_3+a_5+a_6 = 9$.
For maximal weight there is a tie between the first and the last, and either
reduction (i.e.\ using (A) or (D) in Proposition \ref{MN7 basic result}) is
possible.  Note that it is also possible to reduce using (C) (but not (B)), but
the algorithm that we will use in this paper restricts to facets of maximal
weight, so we do not use this option.  The following follows all possible
reductions using facets of maximal weight:  {\scriptsize
\[
\begin{array}{cccccccccccccccccccccccc}
  (3,3,3,1,2,4) \\ \\
   (A) \swarrow \hbox{\hskip 1cm} \searrow (D) \\  \\
(2,2,2,1,2,4) \hbox{\hskip 1.3cm} (3,3,2,1,1,3) \\ \\
(D) \searrow \hbox{\hskip 1cm} \swarrow (A) \\ \\
(2,2,1,1,1,3) \\ \\
\phantom{(C)} \downarrow (C) \\ \\
(2,1,1,0,1,2) \\ \\
   (A) \swarrow \hbox{\hskip 1cm} \searrow (D) \\  \\
(1,0,0,0,1,2) \hbox{\hskip 1.3cm} (2,1,0,0,0,1) \\ \\
(D) \searrow \hbox{\hskip 1cm} \swarrow (A) \\ \\
(1,0,0,0,0,1)
\end{array}
\]
}
\end{example}
\bigskip

We now begin the study of which tetrahedral curves are componentwise
linear.

\begin{theorem} \label{theorem B}
Let $J = (a_1,a_2,a_3,a_4,a_5,a_6)$ be a non-trivial tetrahedral curve.

\begin{itemize}
\item[(a)] If $J$ is minimal then it has a linear resolution, and hence
is componentwise linear.  So from now on we assume that $J$ is not
minimal.

\item[(b)] If (up to permutation of the variables) $J = (0,r,r,r,r,0)$
(i.e. if $J$ has equal non-trivial weights and is supported on a complete
intersection of type $(2,2)$) then $J$ has a pure resolution
that is not linear, and hence is not componentwise linear.

\item[(c)] Suppose that $J$ reduces to another tetrahedral ideal $I$
following the algorithm of \cite{MN7}, i.e.\ using a facet of maximal
weight, and using one of the reductions (A), (B), (C) or (D) of
Proposition \ref{MN7 basic result}.  (If $J$ is not minimal, it is always
possible to reduce using a facet of maximal weight, thanks to Lemma
\ref{reduce max wt lemma}.)  Then the polynomial $F$ that is prescribed by
that algorithm is  a minimal generator of $I$ if and only if $J$ is of the
form described in (b).
\end{itemize}
\end{theorem}

\begin{proof}
Statement (a) follows from \cite{MN7} Theorem 4.2, which in particular
shows that $J$ has a linear resolution.  For (b), we have that $J$ is the
$r$-th power of a complete intersection $(A,B) = (ab, cd)$ of type
(2,2). In fact, it is easy to see that $(ab, cd)^r \subset J$, hence
we get equality because both ideals have the same degree.
Now the
result follows from the fact that its Hilbert-Burch matrix is
\[
\left [
\begin{array}{ccccccccccccccc}
A & B & 0 & 0 & 0 &  & 0 & 0 \\
0 & A & B & 0 & 0 &  & 0 & 0 \\
0 & 0 & A & B & 0 &  & 0 & 0 \\
  &   &   & \vdots \\
0 & 0 & 0 & 0 & 0 &  & A & B

\end{array}
\right ]
\]
where there are $r$ rows and $r+1$ columns, and all entries have degree 2.

For (c), without loss of generality assume that $a_1,a_2,a_3$ give the
facet of maximal weight, so that we use the reduction (A).  We then have
$I = (a_1', a_2', a_3', a_4,a_5,a_6)$ and $F = b^{a_1}c^{a_2}d^{a_3}$.
Suppose first that $F$ is a minimal generator of $I$.  We will show that
then it must be of the type described in (b).

The fact that $F$ is a minimal generator of $I$ means that if we reduce
any of the exponents of $F$, the result is no longer in $I$.  Since $F =
b^{a_1}c^{a_2}d^{a_3}$, it is clear that $F$ vanishes on the components
$(a,b)^{a_1'}$, $(a,c)^{a_2'}$ and $(a,d)^{a_3'}$.  The condition that
$F$ vanishes on $(b,c)^{a_4}$ is given by the inequality $a_1+a_2 \geq
a_4$.  Similarly, the condition that $F$ vanishes on $(b,d)^{a_5}$ is
given by the inequality $a_1+a_3 \geq a_5$ and the condition that $F$
vanishes on $(c,d)^{a_6}$ is given by the inequality $a_2+a_3 \geq a_6$.
To say that reducing any one of the exponents of $F$ by one makes the
result  no longer be in $I$ means that {\em two} of these inequalities
must in fact be equalities. Indeed, this is easily seen if $a_1, a_2, a_3$
are all positive. Assume that without loss of generality $a_1 = 0$ and
only one of these inequalities is an equality. Then, this must be $a_2
+ a_3 = a_6$. And by the assumption for the time being we have $a_2 >
a_4$ and $a_3 > a_5$. But $a_1 + a_2 +a_3$ is the largest weight of a
facet, thus in particular, $a_6 = a_2 + a_3 \geq a_3 + a_5 + a_6$,
hence we get $a_3 = a_5 = 0$, a contradiction.

Therefore, we may assume without loss of generality  that
\[
J = (a_1,a_2,a_3, a_1+a_2, a_1+a_3, a_6).
\]
We assumed that $a_1,a_2,a_3$ give the facet of maximal weight for $J$.
This means, in particular, that $a_1,a_4,a_5$ do not give a facet of
greater weight, i.e.
\[
a_1+(a_1+a_2) + (a_1+a_3) \leq a_1+a_2+a_3.
\]
This forces $a_1 = 0$, and $I = (0, a_2,a_3, a_2,a_3,a_6)$.  Similarly we
have that $a_2,a_4,a_6$ do not give a facet of greater weight for $J$, and
$a_3,a_5,a_6$ do not give a facet of greater weight, so
\[
\left.
\begin{array}{rcl}
a_2+a_2+a_6 & \leq & a_2+a_3 \\
a_3+a_3+a_6 & \leq & a_2+a_3
\end{array}
\right \}
\Rightarrow
\begin{array}{rcl}
a_2+a_6 & \leq & a_3 \\
a_3 + a_6 & \leq & a_2
\end{array}
\]
which means
\[
a_2+2a_6 \leq a_3+a_6 \leq a_2,
\]
so also $a_6 = 0$.  Hence in fact $J$ is of the form $(0,
a_2,a_3,a_2,a_3,0)$.  But then the two inequalities above give $a_2 \leq
a_3$ and $a_3 \leq a_2$, which means that $a_2 = a_3$.  So we have shown
that if we do the reduction via Proposition 3.1 of \cite{MN7} and if
the resulting $F$ is a minimal generator of $I$, then
the ideal $J$ that we started with must be of the type described in (b),
i.e. must be supported on a complete intersection of type (2,2), with
equal weights on each component.  (In particular, $I$ and $J$ must be
arithmetically Cohen-Macaulay.)

Conversely, suppose that $J$ is reduced to $I$ and $J$ is of the form
described in (b).  Without loss of generality say that $J =
(0,r,r,r,r,0)$.  Without loss of generality suppose that we are using
reduction (A).  Then $I$ is  $(0,r-1,r-1,r,r,0)$ and consequently $F$ is
$c^{r}d^{r}$. It is clear that $F$ is a minimal generator of $I$, since
$c^{r-1}d^r$ is not in $(b,c)^r$ and $c^rd^{r-1}$ is not in $(b,d)^r$.
\end{proof}

\begin{corollary} \label{cwl conclusion}
Let $J = (a_1,\dots,a_6)$ be a non-trivial tetrahedral curve.
\begin{itemize}
\item[(a)] If $J$ is not arithmetically Cohen-Macaulay then $J$ is
componentwise linear.

\item[(b)]  Assume that $J$ is arithmetically Cohen-Macaulay.  We can
reduce $J$ to the trivial curve in a finite sequence of steps, each time
using a facet of maximal weight and applying \cite{MN7} Proposition 3.1.
Then the following are equivalent:

\begin{itemize}
\item[(i)] $J$ is componentwise linear;

\item[(ii)] this sequence of steps does not include any curve of the
type described in Theorem
\ref{theorem B} (b);

\item[(iii)]  this sequence of steps
does not include a complete intersection of type (2,2), i.e.  does not
include  any of the curves $(0,1,1,1,1,0)$, $(1,0,1,1,0,1)$, or
$(1,1,0,0,1,1)$.
\end{itemize}
\end{itemize}
\end{corollary}

\begin{proof}
Assume that $J$ is not arithmetically Cohen-Macaulay.  If it is minimal
then by Theorem \ref{theorem B} (a) it is componentwise linear.  If it is
not minimal then we can reduce via facets of maximal weight to a minimal
curve.  In each step, the polynomial $F$ used is not a minimal generator
of the smaller curve $I$, thanks to Theorem \ref{theorem B} (c) and the
fact that $J$ is not arithmetically Cohen-Macaulay.  Then the statement
of (a) follows from Theorem \ref{theorem A} and induction on the number of
steps to a minimal curve.

For (b), the fact that (ii) implies (iii) is trivial, and the implication
(iii) $\Rightarrow$ (ii) follows since a curve of type $(0,r,r,r,r,0)$
reduces to one of type $(0,1,1,1,1,0)$.

Assume that (ii) holds.  Then from Theorem \ref{theorem B} (c), each step
of the procedure of reducing by maximal facets involves a polynomial $F$
that is not a minimal generator of the smaller curve.  Hence by Corollary
\ref{first cor}, each $J$ is componentwise linear if and only if the next
curve $I$ is componentwise linear.  But one can easily check that in
reducing to the trivial curve via facets of maximal weight, eventually one
passes through a curve consisting of all 0's and 1's.  By hypothesis we
do not pass through a complete intersection of type (2,2) (i.e. the curve
$(0,1,1,1,1,0)$, up to permutation).  One can easily check that all other
tetrahedral curves with only entries that are 0 or 1 are componentwise
linear.  Hence by induction the tetrahedral curve $J$ that we started
with is componentwise linear.

Conversely, assume that $J$ is arithmetically Cohen-Macaulay
and componentwise linear.  Again we reduce by facets of maximal weight down
to the trivial curve.  Suppose that at some step we reach a curve of the
type described in Theorem \ref{theorem B} (b), and consider the first
such instance.  We have seen in Theorem \ref{theorem B} (c) that the form
$F$ used in the reduction is a minimal generator if and only if the
larger curve (corresponding to $J$) is of the form described in Theorem
\ref{theorem B} (b).  In our situation we have arrived at the first such
curve, so the larger curve in each step has not been of this form.  Hence
each step in this process has used a form $F$ that was not a minimal
generator of the smaller curve.  Hence by Corollary \ref{first cor},
since we started with an ideal $J$ that was componentwise linear, each of
the smaller curves had ideals $I$ that are also componentwise linear.
But reaching a curve of the type in Theorem \ref{theorem B} (b) we have
obtained one that is not componentwise linear.  This contradiction
completes the proof.
\end{proof}

\begin{definition}
A {\em Schwartau curve} is a tetrahedral curve  $C
= (a_1,a_2,a_3,a_4,a_5,a_6)$ for which $a_2 = a_5 = 0$.
\end{definition}

Note that a Schwartau curve is supported on a complete intersection of
type (2,2).  These curves were studied by P.\ Schwartau in his thesis
\cite{schwartau}.

\begin{corollary} \label{schwartau curves}
Let $J$ be the ideal of a Schwartau
curve.  Then $J$ fails to be componentwise linear if
and only if all of $a_1,a_3,a_4,a_6$ are $>0$ and $a_1+a_6 = a_3+a_4$.
\end{corollary}

\begin{proof}
We reduce $J$ to the trivial curve by a sequence of steps using the
reduction of \cite{MN7}, Proposition 3.1, and using facets of maximal
weight.  By Corollary \ref{cwl conclusion}, $J$ fails to be componentwise
linear if and only if this reduction includes the curve $(1,0,1,1,0,1)$
(this time it must be precisely this curve, not up to permutation).

Suppose that $J$ fails to be componentwise linear.  If any of the $a_i$
are 0, then clearly we cannot hope to reach $(1,0,1,1,0,1)$.  But note
that each step in the reduction reduces both sums $a_1+a_6$ and $a_3+a_4$
by 1, so if these sums are not equal to begin with, they will never be
equal.  Hence we will never reach $(1,0,1,1,0,1)$.  Hence we must have
the claimed equality.

Conversely, assume that all $a_i > 0$ and that $a_1+a_6 = a_3+a_4$.  The
maximal facet will always include $\max \{ a_1,a_6 \}$ and $\max \{
a_3,a_4 \}$ (and the third edge is 0).  Hence since $a_1+a_6 = a_3+a_4$,
we eventually arrive at $(1,0,1,1,0,1)$, so $J$ is not componentwise
linear.
\end{proof}

Note that the curves considered in Corollary \ref{schwartau curves} are
automatically arithmetically Cohen-Macaulay, thanks to Corollary \ref{cwl
conclusion} (a).

\begin{corollary} \label{hope}
Let $(a_1,a_2,a_3,a_4,a_5,a_6)$ be a tetrahedral curve $C$.  Consider the
sums  $a_1+a_6, a_2+a_5,
a_3+a_4$.  If  the curve  fails to be componentwise linear,
then the two larger of these sums are equal.
\end{corollary}

\begin{proof}
We know that if $C$ is not componentwise linear, then it is
arithmetically Cohen-Macaulay and reduces via
facets of maximal weight to one of the curves listed in Corollary
\ref{cwl conclusion} (b)(iii). Notice that for any of these curves the
two larger sums equal two and the third is zero. But each basic double
link increases the two larger sums by one and the third sum by zero or
one. Hence, the claim for $C$ follows.
\end{proof}

The converse to this statement is not true. Here is a counterexample.

\begin{example}
Let $J = (10, 1, 2, 3, 10, 1)$. Then this curve is arithmetically
Cohen-Macaulay, but is componentwise linear because the reduction to
the trivial curve does not pass through any of the curves listed in
Corollary \ref{cwl conclusion} (b).
\end{example}

%%%%%%%%%%%%%%%%%%%%%%%%%%%%%%%%%%%%%%%%%%%%%%%%%%%%%%%%%%%%%%%%%%%%%%

\section{The Minimal Free Resolution of a Tetrahedral Curve}

In this section we describe the whole minimal free
resolution of a tetrahedral curve.  In particular, we make observations about
the minimal generators and about the regularity.

We first prove a lemma about the degree of the monomial $F$ used in the
algorithm from \cite{MN7} for reducing a tetrahedral curve. We focus on the
case in which our ideals are not componentwise linear.

\begin{lemma} \label{Fdegree} Let $J$ be the ideal of a tetrahedral curve
$(a_1,\dots,a_6)$ that is not of type $(0,r,r,r,r,0)$, even with the variables
permuted. Assume that $J$ reduces using the algorithm
from \cite{MN7} to a curve of type $(0,r,r,r,r,0)$ (possibly with the variables
permuted), and hence is not componentwise linear.  Assume further that $J$ has
its lowest degree minimal generators in degree $p$. If $I$ is the ideal
of the curve obtained by reducing the facet of maximal weight of the curve of
$J$ using the algorithm from
\cite{MN7}, and
$J=L
\cdot I + (F)$, then $\deg F \ge p+1$.
\end{lemma}

\begin{proof} Note that $I$ has its lowest degree minimal generators in degree
$p-1$ because $F$ is not a minimal generator of $I$. Therefore $F$ 
has degree at
least $p$, and we wish to show that it has degree at least $p+1$. 
Suppose to the
contrary that $\deg F=p$ and that $F=b^{a_1}c^{a_2}d^{a_3}$. Then $F$ is $b$,
$c$, or $d$ times a minimal generator of $I$; without loss of generality, say
$b^{a_1-1}c^{a_2}d^{a_3}$ is a minimal generator of $I$. Then
$b^{a_1-1}c^{a_2-1}d^{a_3} \not \in I$, so $a_1-1+a_2-1 < a_4$ or $a_2-1+a_3 <
a_6$. Similarly, $b^{a_1-1}c^{a_2}d^{a_3-1} \not \in I$, and thus 
$a_1-1+a_3-1 <
a_5$ or $a_2-1+a_3 < a_6$.

Suppose first that $a_2+a_3-1 < a_6$; then $a_2+a_3=a_6$ because
$b^{a_1-1}c^{a_2}d^{a_3} \in I$. We are assuming that $a_1+a_2+a_3$ is the
maximal weight of a facet, and therefore $a_1+a_6=a_1+a_2+a_3 \ge a_2+a_4+a_6$.
Consequently, $a_1 \ge a_2+a_4$. Also, we have $a_1 \ge a_4 - a_2$ because
$b^{a_1}c^{a_2}d^{a_3} \in I$, and hence $a_1 \ge a_4$. Similarly, using that
$a_1+a_2+a_3 \ge a_3+a_5+a_6$, we conclude that $a_1 \ge a_5$. As a result,
$$
a_1 + a_6 = a_1+a_2+a_3 \ge a_5+a_2+a_3 \quad \hbox{and} \quad a_1 + a_6 \ge
a_4+a_2+a_3.
$$
This says that $a_1+a_6$ is equal to the maximum of
$\{a_1+a_6,a_2+a_5,a_3+a_4\}$.

Because $J$ reduces to a curve of type $(0,r,r,r,r,0)$, it is not componentwise
linear. Hence by Corollary \ref{hope}, $a_1+a_6$ is equal to either 
$a_2+a_5$ or
$a_3+a_4$; without loss of generality, assume it is $a_2+a_5$. We have
$$
a_1+a_2+a_3 = a_1+a_6 =a_2 + a_5,
$$
so $a_1+a_3=a_5$. Since $a_1 \ge a_5$, this forces $a_3=0$. Therefore
$a_2=a_6$, and because $a_1+a_6=a_2+a_5$, $a_1=a_5$. Hence $J$ is the 
ideal of a
curve $(a_1,a_2,0,a_4,a_1,a_2)$. But then $I$ is the ideal of a curve 
of the form
$(a_1-1,a_2',0,a_4,a_1,a_2)$. We know that $b^{a_1-1}c^{a_2}d^{a_3}=$
$b^{a_1-1}c^{a_2}$ is a minimal generator of $I$, but
$$
b^{a_1-1}c^{a_2} \not \in (b,d)^{a_1},
$$
a contradiction.

As a result, we conclude that $a_1+a_3-1=a_5$ and $a_1+a_2-1=a_4$. Therefore
$$
a_2+a_5 = (a_4-a_1+1)+a_5=(a_4-a_1+1)+(a_1+a_3-1)=a_3+a_4.
$$
Because $J$ is not componentwise linear, $a_2+a_5=a_3+a_4$ must be the maximum
value among $a_1+a_6$, $a_2+a_5$, and $a_3+a_4$ since by Corollary \ref{hope},
the largest two of those are equal. Using this and the fact that 
$a_1+a_3=a_5+1$,
we have
$$
a_3+a_4+1=a_2+a_5+1=a_1+a_2+a_3 \ge a_1+a_4+a_5,
$$
where the inequality holds because $a_1+a_2+a_3$ gives the maximal weight of a
facet. Hence
$$
a_3+1 \ge a_1+a_5 = a_1 + (a_1+a_3-1),
$$
so $2 \ge 2a_1$, and $a_1 \le 1.$ But $a_1 \not = 0$ because if it were zero,
$b^{a_1-1}c^{a_2}d^{a_3}$ would not be a minimal generator of $I$. Therefore
$a_1=1$, which implies that $a_3=a_5$ and $a_2=a_4$. Thus $J$ is the ideal of a
curve $(1,a_2,a_3,a_2,a_3,a_6)$, and $F=bc^{a_2}d^{a_3}$. Note that 
$a_2 \not = 0
\not = a_3$, for if one of them were zero, then $J$ could not reduce to a curve
of the form $(0,r,r,r,r,0)$, even with the variables permuted. 
Consequently, $I$
is the ideal of a curve $(0,a_2-1,a_3-1,a_2,a_3,a_6)$.

Next, we claim that $a_6 \le 1$. To see this, note that $1+a_2+a_3$ gives the
maximal weight of a facet of $J$. Therefore
$$
1+a_2+a_3 \ge 2a_2+a_6 \quad \hbox{and} \quad 1+a_2+a_3 \ge 2a_3+a_6;
$$
adding these inequalities gives the claim.

We wish to show that $m=bc^{a_2-1}d^{a_3-1} \in I$. If so, then it 
has degree at
least $p-1$, which implies that $\deg F = \deg bc^{a_2}d^{a_3} \ge p+1$. If
$a_2+a_3-2 \ge a_6$, we can conclude that $m \in I$. This inequality holds if
$a_6=0$ since $a_2$ and $a_3$ are both at least one. If $a_2+a_3-2 < 
a_6=1$, then
$a_2=a_3=1$, and $I$ is the ideal of the curve $(0,0,0,1,1,1)$, which does not
reduce to an ideal of the form $(0,r,r,r,r,0)$. Consequently, $\deg F \ge p+1$.
\end{proof}

The lemma allows us to compare the resolutions of $J$ and $\gin(J)$ 
for any ideal
$J$ of a tetrahedral curve.

\begin{proposition} \label{ginres} Let $J \subset k[a,b,c,d]$ be the ideal of a nontrivial tetrahedral curve, and suppose the characteristic of $k$ is zero.

\begin{itemize}

\item[(a)] If $J$ is componentwise linear, then the graded Betti numbers of $J$
and $\gin(J)$ are the same.

\item[(b)] Suppose $J$ is not componentwise linear and has its lowest degree
minimal generators in degree $p$. Assume that $J$ reduces using the algorithm
from \cite{MN7} to a curve of type $(0,r,r,r,r,0)$ (possibly with the variables
permuted), with $r>0$, but not $(0,r+1,r+1,r+1,r+1,0)$. Then the graded Betti
numbers of $\gin(J)$ and $J$ are the same except that $\gin(J)$ has $r$
additional minimal generators and syzygies in degree $p+1$.
\end{itemize}
\end{proposition}

\begin{proof} Part (a) is immediate from Theorem 1.1 in \cite{AAH}. 
For part (b),
note that $J$ has projective dimension two because $R/J$ is Cohen-Macaulay.
Suppose first that $J$ is the ideal of a $(0,r,r,r,r,0)$ curve. Then it is easy
to compute (see Theorem \ref{theorem B}(b)) that $J$ has resolution
$$
0 \rightarrow R(-2r-2)^r \rightarrow R(-2r)^{r+1} \rightarrow J \rightarrow 0.
$$
Since the regularity, Hilbert functions, and projective dimensions of $J$ and
$\gin(J)$ are the same, the only possible differences in their Betti 
numbers are
additional generators and syzygies of $\gin(J)$ in degree $2r+1$. Because
$\gin(J)$ is strongly stable, it is componentwise linear, and thus
$(\gin(J)_{2r})$ has a linear resolution. Therefore it must have $r$ minimal
syzygies of degree $2r+1$ on the $r+1$ minimal generators of degree $2r$. Thus
there are $r$ additional generators of degree $2r+1$ to preserve the Hilbert
function. Note that $(J_{2r+1+s})$ has a linear resolution for all $s \ge 0$
since the regularity of $J$ is $2r+1$.

Suppose now that $J$, the ideal of a curve $(a_1,\dots,a_6)$, is a basic double
link of $I$, so that $J=L \cdot I+(F)$, where $L$ is a linear form. We assume
that $I$ is obtained from $J$ by reducing a facet of maximal weight. Suppose
further that $J$ is not a curve of type $(0,r,r,r,r,0)$ but reduces 
to a curve of
that form (again possibly with the variables permuted) and not
$(0,r+1,r+1,r+1,r+1,0)$. We have that $I$ has its lowest degree minimal
generators in degree $p-1$ because $F$ is not a minimal generator of $I$.  For
the induction hypothesis, we assume that $I$ has $r+1$ minimal generators of
degree $p-1$, no minimal syzygies of degree $p$, and that $(I_{p+s})$ has a
linear resolution for all $s \ge 0$.

By Lemma \ref{Fdegree}, $\deg F \ge p+1$. Therefore $J$ has $r+1$ minimal
generators of lowest degree $p$ and no syzygies of degree $p+1$. Because $F$ is
not a minimal generator of $I$, it follows from the same arguments as 
in Theorem
\ref{theorem A} and Corollary \ref{where resol is not lin} that since 
$(I_{p+s})$
has a linear resolution for all $s \ge 0$, $(J_{p+1+s})$ does also. That means
that $J_{\ge (p+1)}$ is componentwise linear, where $J_{\ge (p+1)}$ 
is the ideal
generated by all elements of $J$ with degree at least $p+1$.

By Lemma \ref{truncatebetti}, $\beta_{i,i+j}(J_{\ge (p+1)}) = \beta_{i,i+j}(J)$
for all $j \ge p+2$. Consequently, $\beta_{i,i+j}(J) = \beta_{i,i+j}(\gin(J))$
for all $j \ge p+2$ (and for all $j < p$ since those Betti numbers are zero).
Because the gin preserves the Hilbert function, we have $\beta_{0,p}(J) =
\beta_{0,p}(\gin(J))$, and also $\beta_{1,p+2}(J) = 
\beta_{1,p+2}(\gin(J))$ since
the number of generators of degree $p+2$ is the same for both ideals. As a
result, the only possible changes are in degree $p+1$. But $\gin(J)$ 
is strongly
stable, and thus since $J$ has $r+1$ minimal generators of degree 
$p$, $\gin(J)$
must have $r$ syzygies of degree $p+1$. Consequently, $\gin(J)$ has $r$
additional minimal generators of degree $p+1$ to preserve the Hilbert function.
\end{proof}

\begin{example} \label{ginprocess} Let $J$ be the ideal of the 
tetrahedral curve
$(2,5,5,5,5,0)$. Then $J$ reduces to the curve $(0,4,4,4,4,0)$, and $J$ has
resolution
$$
0 \rightarrow R(-13)^2 \oplus R(-12)^4 \rightarrow R(-12)^2 \oplus R(-10)^5
\rightarrow J \rightarrow 0.
$$

The gin of $J$ must add four minimal generators and syzygies of 
degree 11, and it
has resolution
$$
0 \rightarrow R(-13)^2 \oplus R(-12)^4 \oplus R(-11)^4 \rightarrow R(-12)^2
\oplus R(-11)^4 \oplus R(-10)^5 \rightarrow \gin(J) \rightarrow 0.
$$
\end{example}

As a consequence, we get a description of the regularity of the ideal of any
tetrahedral curve in terms of the degree of the highest degree 
minimal generator.

\begin{corollary} \label{regularity} Let $J$ be the ideal of a nontrivial
tetrahedral curve.

\begin{itemize}

\item[(a)] If $J$ is the ideal of a curve of the form $(0,r,r,r,r,0)$, then the
regularity of $J$ is $2r+1$.

\item[(b)] If $J$ is not the ideal of a curve of the form $(0,r,r,r,r,0)$
(possibly with the variables permuted), then the regularity of $J$ is 
the degree
of the largest degree minimal generator of $J$.
\end{itemize}
\end{corollary}

\begin{proof} Part (a) is immediate from the resolution of an ideal 
of a curve of
that form. Part (b) is clear when $J$ is componentwise linear, so suppose that
$J$ is not componentwise linear. Note that if $p$ is the smallest 
degree in which
$J$ has minimal generators, then $J$ also has generators in a higher 
degree $p+s$
by Lemma \ref{Fdegree}. Because $J_{\ge p+1}$ is componentwise linear, the
regularity of $J_{\ge p+1}$ is equal to the highest degree in which it has a
minimal generator, and these invariants are the same as for $J$.
\end{proof}

We can use Corollary \ref{regularity} to get a more precise statement about the
regularity of a tetrahedral curve that can be read directly from the 6-tuple
$(a_1,\dots,a_6)$. First, we prove a lemma that will serve as the 
inductive step
in our next result.

\begin{lemma} \label{maxdegree-gen} Let $J$ be the ideal of a nonminimal
tetrahedral curve. Suppose $J=L \cdot I + (F)$ is a basic double link of $I$,
where $I$ is obtained from $J$ by reducing a facet of maximal weight. 
Assume also that the maximal degree of a minimal generator of $I$ is equal to the maximal weight of a facet of $I$. Then $\deg F$ is the highest degree in which $J$ has a minimal generator, and this degree is equal to the maximal weight of a facet of $J$. 

\begin{proof} Note first that the ideals corresponding to $(0,r,r,r,r,0)$ (possibly with the variables permuted) have
all their minimal generators in degree $2r$, and thus the result 
holds if $J$ is
of that form. So suppose that $J$ corresponds to the curve
$(a_1,a_2,a_3,a_4,a_5,a_6)$, which is not of the form $(0,r,r,r,r,0)$, and
$a_1+a_2+a_3$ gives the maximal weight of a facet of this curve. 
Then the curve
corresponding to $I$ has the form $(a_1',a_2',a_3',a_4,a_5,a_6)$, where
$a_i'=\max\{0,a_i-1\}$. We have $J=L \cdot I + (F)$, where $\deg 
F=a_1+a_2+a_3$.
Thus we need to show that $a_1+a_2+a_3$ is at least as large as the maximal
weight of a facet of the curve corresponding to $I$ plus one (adding 
one because
$L$ is a linear form).

If the maximal weight of a facet of $I$ is $a_1'+a_2'+a_3'$, then it is clear
that $a_1+a_2+a_3 \ge a_1'+a_2'+a_3'+1$, so $F$ is a minimal generator of
$J$ of highest degree. Suppose the maximal weight of a facet for $I$ is
$a_1'+a_4+a_5$; the other two remaining cases are the same. Suppose 
$a_1+a_2+a_3
< a_1'+a_4+a_5+1$. Then
$$
a_1+a_2+a_3 \le a_1'+a_4+a_5 \le a_1+a_4+a_5.
$$
Because $a_1+a_2+a_3 \ge a_1+a_4+a_5$, all the inequalities are equalities,
which means that $a_1'=a_1=0$ and $a_2+a_3=a_4+a_5$. Therefore 
$a_1+a_4+a_5$ also
gives the maximal weight of a facet of $J$, and hence
$$
a_4+a_5 \ge a_2 + a_4 + a_6 \quad \hbox{and} \quad a_4+a_5 \ge a_3+a_5+a_6.
$$
Adding these inequalities together and using the fact that $a_2+a_3=a_4+a_5$,
we have $2a_6 \le 0$, so $a_6=0$. Thus $J$ has the form 
$(0,a_2,a_3,a_4,a_5,0)$.

Because $0+a_2+a_3$ gives the maximal weight of a facet of $J$, we have
$$
a_2+a_3 \ge a_3+a_5 \quad \hbox{and} \quad a_2+a_3 \ge a_2 + a_4,
$$
so $a_2 \ge a_5$ and $a_3 \ge a_4$. But $a_2+a_3=a_4+a_5$, so $a_2=a_5$ and
$a_3=a_4$, and $J$ has the form $(0,a_2,a_3,a_3,a_2,0)$.

If $a_2=a_3$, then $J$ is of the form $(0,r,r,r,r,0)$, contradicting our
assumption that it was not. Otherwise, by Corollary \ref{s-min}, $J$ 
is minimal,
which is again a contradiction.
\end{proof}
\end{lemma}

We can now express the regularity of any tetrahedral curve explicitly.

\begin{theorem} \label{reg-weight} Let $J$ be the ideal of a nontrivial
tetrahedral curve.

\begin{itemize}

\item[(a)] If $J$ is the ideal of a curve of the form $(0,r,r,r,r,0)$, possibly
with the variables permuted, then the regularity of $J$ is $2r+1$.

\item[(b)] If $J$ is the ideal of a minimal curve $(a_1,\dots,a_6)$, assume
without loss of generality that $a_6$ is the largest of the $a_i$. Then the
regularity of $J$ is $a_1+a_6$, which is strictly greater than the weight of a
maximal facet.

\item[(c)] Suppose $J$ is the ideal of a curve $(a_1,\dots,a_6)$ that is not
minimal and  not of the form in (a). Then the regularity of $J$ is the maximal
weight of a facet.
\end{itemize}
\end{theorem}

\begin{proof} Part (a) follows from computing the resolution of ideals of this
type; see Theorem \ref{theorem B}(b). Consider now part (b). The value of the
regularity is an immediate consequence of Theorem \ref{thm-res}. For the
inequality, assume without loss of generality that $a_1+a_2+a_3$ gives the
maximal weight of a facet; the argument is similar if $a_6$ is included in the
maximal weight of a facet. By Corollary \ref{s-min}, $a_6 > a_2+a_3$, and thus
$a_1+a_6 > a_1+a_2+a_3$.

We turn now to part (c). First, we consider the case in which $J$ is
arithmetically Cohen-Macaulay but not of the form $(0,r,r,r,r,0)$, 
even with the
variables permuted. It is easy to check that all curves of weight at most three
and all $a_i$ zero or one have regularity equal to the maximal weight 
of a facet;
also, the regularity is equal to the highest degree of a minimal generator in
each case.  Moreover, all arithmetically Cohen-Macaulay $J$ reduce to these
cases. We proceed by induction with these curves as the base case.

We reduce $J$ down to a base case with the algorithm from \cite{MN7}, always
reducing by a facet of maximal weight. For the induction hypothesis, assume the
following: All ideals $M$ below $J$ in the reduction from $J$ down to 
a base case
have the property that the highest degree of a minimal generator of 
$M$ is equal
to the maximal weight of a facet of $M$. Suppose that $I$ is obtained 
by reducing
a facet of $J$ of maximal weight. By the induction hypothesis, the 
maximal degree
of a minimal generator of $I$ is the maximal weight of a facet of $I$. By Lemma
\ref{maxdegree-gen}, this implies that the maximal degree of a 
minimal generator
of $J$ is the maximal weight of a facet of $J$. We conclude from Corollary
\ref{regularity} that this quantity is equal to the regularity of $J$.

Finally, we consider the case in which $J$ is not arithmetically 
Cohen-Macaulay.
Suppose first that $I$ is the ideal of a minimal, not arithmetically
Cohen-Macaulay curve, and suppose $I$ is obtained from $J$ by 
reducing a facet of
maximal weight as in the algorithm from \cite{MN7}. Then $J=L \cdot I + (F)$,
where $L$ is a linear form. Without loss of generality, suppose that $J$ is the
ideal of a curve $(a_1,\dots,a_6)$ with $\deg F=a_1+a_2+a_3$ giving the maximal
weight of a facet. Then $I$ is the ideal of a curve
$(a_1',a_2',a_3',a_4,a_5,a_6)$. We may assume that
$a_6=\max\{a_1',a_2',a_3',a_4,a_5,a_6\}$; the argument is similar if 
the maximal
weight of a facet of $I$ includes $a_6$.

The regularity of $I$ is $a_1'+a_6$, and thus the minimal generators 
of $L \cdot
I$ have degree $a_1'+a_6+1$. Hence $\text{reg}(J) \ge a_1'+a_6+1$. By Corollary
\ref{regularity}, the regularity of $J$ is the maximal degree of a minimal
generator. Therefore
$$
\text{reg}(J)=\max\{\deg F=a_1+a_2+a_3,a_1'+a_6+1\}.
$$
We want to show that $a_1+a_2+a_3 \ge a_1'+a_6+1$.

Initially, note that $a_1'=a_1-1$; otherwise, $a_1'=0$, contradicting Corollary
\ref{s-min}. Thus we need to show that $a_2+a_3 \ge a_6$. Suppose instead that
$a_2+a_3 < a_6$. We show that $J$ is then the ideal of a minimal 
curve, which is
a contradiction.

Because $I$ is the ideal of a minimal curve, $a_6 > a_4+a_5$. Thus
$$
a_6 > \max\{a_4+a_5, a_2+a_3\}.
$$
Now, since $a_1+a_2+a_3$ gives the maximal weight of a facet of $J$, we have
$a_1+a_2+a_3 \ge a_2+a_4+a_6$. Consequently,
$$
a_1+a_2+a_3 \ge a_2+a_4+a_6 > a_2+a_4+(a_2+a_3),
$$
using the assumption that $a_2+a_3>a_6$. Therefore $a_1 > a_2+a_4$. Similarly,
$$
a_1+a_2+a_3 \ge a_3+a_5+a_6 > a_3+a_5+(a_2+a_3),
$$
so $a_1 > a_3+a_5$. Hence
$$
a_1 > \max\{a_3+a_5, a_2+a_4\}.
$$
Combining this with the inequalities for $a_6$ and using Corollary
\ref{s-min}, we conclude that $J$ is the ideal of a minimal curve, a
contradiction.

Thus if $J$ reduces in one step to the ideal of a minimal curve by reducing a
facet of maximal weight, the maximal degree of a minimal generator of $J$ is
equal to the maximal weight of a facet of $J$. By Corollary \ref{regularity},
this is also the regularity of $J$. Now the result for all ideals of 
nonminimal,
non-arithmetically Cohen-Macaulay curves follows from the same 
induction process
as in the arithmetically Cohen-Macaulay case.
\end{proof}

With Theorem \ref{reg-weight}, we can read the regularity of a 
tetrahedral curve
straight from the $a_i$ that describe it. With a bit more work, we 
can also find
the graded Betti numbers of the ideal of a tetrahedral curve without any
substantial computation. We begin with a lemma that describes how the maximal
weight of a facet changes when we do a basic double link. Its proof 
is similar to
several of our earlier arguments in this section.

\begin{lemma} \label{maxweight-increases} Let $J$ be the ideal of a nonminimal
tetrahedral curve that is not of the form $(0,r,r,r,r,0)$, even with the
variables permuted. Suppose $J$ is a basic double link of $I$, and $I$ is
obtained by reducing a facet of maximal weight. Then the maximal weight of a
facet of
$J$ is strictly larger than the maximal weight of a facet of $I$.
\end{lemma}

\begin{proof} Let $J$ be the ideal of the tetrahedral curve
$(a_1,a_2,a_3,a_4,a_5,a_6)$, and assume without loss of generality that
$a_1+a_2+a_3$ gives the maximal weight of a facet. Then $I$ is the ideal of a
curve $(a_1',a_2',a_3',a_4,a_5,a_6)$, where $a_i'=\max\{0,a_i-1\}$. Clearly
$a_1+a_2+a_3$ is at least as large as the weight of any facet of $I$, 
but suppose
it is equal to the weight of a facet of $I$. Without loss of generality, assume
that $a_1+a_2+a_3=a_1'+a_4+a_5$. Since $a_1+a_2+a_3 \ge a_1+a_4+a_5$, 
we conclude
that $a_1=a_1'=0$, and $a_2+a_3=a_4+a_5$. Because $a_1'+a_4+a_5=a_4+a_5$ gives
the maximal weight of a facet of $I$, we have
$$
a_4+a_5 \ge a_2'+a_4+a_6 \quad \hbox{and} \quad a_4+a_5 \ge a_3'+a_5+a_6,
$$
so $a_5 \ge a_2'+a_6$ and $a_4 \ge a_3'+a_6$. Adding these inequalities and
using the fact that $a_2+a_3=a_4+a_5$, we have
$$
a_2+a_3 = a_4+a_5 \ge a_2'+a_3'+2a_6,
$$
and thus $a_6 \le 1$.

Suppose first that $a_6=1$. Because $a_2+a_3$ is the maximal weight 
of a facet of
$J$, we have
$$
a_2+a_3 \ge a_2+a_4+1 \quad \hbox{and} \quad a_2+a_3 \ge a_3+a_5+1,
$$
which implies that $a_3 \ge a_4+1$ and $a_2 \ge a_5+1$. This contradicts the
fact that $a_2+a_3=a_4+a_5$.

Therefore $a_6=0$. Then
$$
a_4+a_5=a_1+a_4+a_5=a_1+a_2+a_3 \ge a_2+a_4+a_6 = a_2+a_4;
$$
hence $a_5 \ge a_2$. Similarly, $a_4 \ge a_3$. Because $a_2+a_3=a_4+a_5$,
$a_5=a_2$ and $a_4=a_3$, and $J$ has the form $(0,a_2,a_3,a_3,a_2,0)$. If
$a_2=a_3$, this contradicts the assumption that $J$ is not of the form
$(0,r,r,r,r,0)$. Otherwise, by Corollary \ref{s-min}, $J$ is the ideal of a
minimal curve, again a contradiction.
\end{proof}

As a corollary, we obtain a characterization of when the mapping cone 
resolution
of $J$ coming from a basic double link is minimal.

\begin{corollary} \label{cone} Let $J$ be the ideal of a nonminimal tetrahedral
curve. Assume that  $J$ is a basic double link of $I$, so $J=L \cdot I + (F)$,
where $L$ is a linear form, and assume that $I$ is obtained by reducing a facet
of
$J$ of maximal weight. Set $\deg F = e$. Then the mapping cone 
resolution of $J$
coming from the short exact sequence
$$
0 \rightarrow R(-e-1) \rightarrow I(-1) \oplus R(-e) \rightarrow J \rightarrow
0
$$
is minimal if and only if $J$ is not the ideal of a curve $(0,r,r,r,r,0)$
(possibly with the variables permuted).
\end{corollary}

\begin{proof} First, suppose $J$ is of the form $(0,r,r,r,r,0)$. Then 
by Theorem
\ref{theorem B}, $F$ is a minimal generator of $I$. Therefore the mapping cone
resolution cannot be minimal, for $L \cdot F$ will not be one of the minimal
generators of $J$.

Now suppose that $J$ is not the ideal of a curve $(0,r,r,r,r,0)$. By Lemma
\ref{maxweight-increases}, the maximal weight of a facet of $J$ is strictly
greater than the maximal weight of a facet of $I$. If $I$ is not the ideal of a
minimal, non-arithmetically Cohen-Macaulay curve, the maximal weight of a facet
of $I$ is equal to the maximal degree in which $I$ has a minimal generator by
Lemma \ref{reg-weight}. Therefore the degree of $F$ is at least as large as the
degree of the highest degree of a minimal generator of $L \cdot I$. Because all
the minimal generators of $L \cdot I$ are divisible by $L$, and $F$ is not, $F$
is not a redundant generator of $J$, and its degree is too high to make any of
the minimal generators of $L \cdot I$ redundant. Hence the mapping cone
resolution is minimal.

If $I$ is not arithmetically Cohen-Macaulay, and $I$ is minimal, then 
$F$ is not
a minimal generator of $I$ by Theorem \ref{theorem B}. Therefore the degree of
$F$ is at least as large as the degree of the minimal generators of 
$L \cdot I$,
and again, $F$ and all of the minimal generators of $L \cdot I$ are minimal
generators of $J$. This proves that the mapping cone resolution is minimal.
\end{proof}

\begin{remark} \label{minimalres} We can use Corollary \ref{cone} to help
describe an inductive procedure with which we can easily compute the 
graded Betti
numbers of the ideal of any tetrahedral curve. Suppose $J$ is the ideal of a
tetrahedral curve. Using the algorithm from \cite{MN7}, reducing by a facet of
maximal weight, we get a sequence of reductions
$$
J=J_s \mapsto J_{s-1} \mapsto \cdots \mapsto J_1 \mapsto J_0=M.
$$
If $J$ is not arithmetically Cohen-Macaulay, let $M$ be the ideal of the
minimal curve to which the curve corresponding to $J$ reduces. If $J$ is
arithmetically Cohen-Macaulay and componentwise linear, let $M$ be the ideal of
the trivial curve; that is, $M=(1)$. Finally, if $J$ is arithmetically
Cohen-Macaulay and not componentwise linear, suppose $J$ reduces to 
the ideal of
a curve of the form $(0,r,r,r,r,0)$ but not $(0,r+1,r+1,r+1,r+1,0)$, 
and let $M$
be the ideal of $(0,r,r,r,r,0)$. In all three cases, we know the minimal graded
free resolution of $M$ (in the nontrivial cases, from Theorem \ref{thm-res} or
Theorem \ref{theorem B}(b)). By Corollary \ref{cone}, the mapping 
cone resolution
of $J_r$ obtained from the short exact sequence induced by the basic 
double link
$J_r = L_r \cdot J_{r-1} + (F_r)$ is minimal. Therefore to get the minimal
resolution of $J_r$, one shifts the minimal resolution of $J_{r-1}$ 
by one degree
and adds a generator of degree $\deg F_{r}$ and a syzygy of degree $(\deg
F_{r}+1)$. If $J \not = M$, we can read the maximal degree of a 
minimal generator
(and the corresponding highest degree first syzygy) directly from the maximal
weight of a facet of $J$ and then continue the process inductively 
with the rest
of the $J_r$.
\end{remark}

Once we know the reduction sequence, the minimal free resolution of $J$ can be
written immediately only from knowledge of the sequence and of $M$.  We
illustrate this process in three examples.

\begin{example} \label{ex:acm-cwl} Suppose $J$ is the ideal of the curve
$(1,2,1,2,0,2)$. We illustrate the reduction procedure and degrees of 
generators
and syzygies at each step.

\bigskip

\begin{tabular}{cccc} {\bf Curve} & {\bf Maximal weight} & {\bf Degree of
generator} & {\bf Degree of syzygy}\\ \hline
$(1,2,1,2,0,2)$ & 6 & 6 & 7\\
$(1,1,1,1,0,1)$ & 3 & 3+1=4 & 4+1=5\\
$(0,0,0,1,0,1)$ & 2 & 2+2=4 & 3+2=5\\
$(0,0,0,0,0,0)$\\
\end{tabular}

\bigskip

We add in the resolution of $R$ itself, shifted in degree by three 
because of the
three reductions. Therefore the minimal graded free resolution of $J$ is
$$
0 \rightarrow R(-7) \oplus R(-5)^2 \rightarrow R(-6) \oplus R(-4)^2 \oplus
R(-3) \rightarrow J \rightarrow 0.
$$
\end{example}

\begin{example} \label{ex:acm-notcwl} Let $J$ be the ideal of the curve
$(1,3,4,2,3,0)$. Then $J$ is arithmetically Cohen-Macaulay and not 
componentwise
linear; thus it reduces to a curve of the form $(0,r,r,r,r,0)$, and we know the
resolutions of those curves by Theorem \ref{theorem B}(b).

\bigskip

\begin{tabular}{cccc} {\bf Curve} & {\bf Maximal weight} & {\bf Degree of
generator} & {\bf Degree of syzygy}\\ \hline
$(1,3,4,2,3,0)$ & 8 & 8 & 9\\
$(0,2,3,2,3,0)$ & 6 & 6+1=7 & 7+1=8\\
$(0,2,2,2,2,0)$ \\
\end{tabular}

\bigskip

We add to this the resolution of $(0,2,2,2,2,0)$, shifted by two 
since there were
two reductions. Using Theorem \ref{theorem B}(b), this gives three 
generators of
degree $4+2=6$ and two syzygies of degree $6+2=8$.  Hence the minimal 
resolution
of
$J$ is
$$
0 \rightarrow R(-9) \oplus R(-8)^3 \rightarrow R(-8) \oplus R(-7) \oplus
R(-6)^3 \rightarrow J \rightarrow 0.
$$
\end{example}

\begin{example} \label{ex:notacm} We consider a curve that is not 
arithmetically
Cohen-Macaulay. Let $J$ be the ideal of the curve $(7,5,5,2,1,6)$.

\bigskip

\begin{tabular}{cccc} {\bf Curve} & {\bf Maximal weight} & {\bf Degree of
generator} & {\bf Degree of first syzygy}\\ \hline
$(7,5,5,2,1,6)$ & 17 & 17 & 18\\
$(6,4,4,2,1,6)$ & 14 & 14+1=15 & 15+1=16\\
$(5,3,3,2,1,6)$ & 11 & 11+2=13 & 12+2=14\\
$(4,2,2,2,1,6)$ & 10 & 10+3=13 & 11+3=14\\
$(4,1,2,1,1,5)$
\end{tabular}

\bigskip

To the generators and syzygies from the reduction, we add the resolution of the
minimal curve $(4,1,2,1,1,5)$, shifted by four because of the four 
reductions. By
Theorem \ref{thm-res}, the resolution of the ideal $I$ of $(4,1,2,1,1,5)$ is
$$
0 \rightarrow R(-11)^{14} \rightarrow R(-10)^{37} \rightarrow R(-9)^{24}
\rightarrow I \rightarrow 0.
$$

Therefore the resolution of $J$ is
$$
0 \rightarrow R(-15)^{14} \rightarrow R(-18) \oplus R(-16) \oplus R(-14)^{39}
\rightarrow R(-17) \oplus R(-15) \oplus R(-13)^{26} \rightarrow J 
\rightarrow 0.
$$

\end{example}

In case of non-arithmetically Cohen-Macaulay curves, part of the
preceding discussion can be summarized as follows.

\begin{corollary} \label{cor-mfr}
Let $J$ be the ideal of a tetrahedral curve that is not arithmetically
Cohen-Macaulay. Then its minimal free resolution is of the form
$$
0 \to R^{\beta_3}(-e_0 - s - 2) \to \begin{array}{c} G(-1) \\
\oplus \\
R^{\beta_2}(-e_0 - s - 1)
\end{array}
\to \begin{array}{c}
G \\
\oplus \\
 R^{\beta_1}(-e_0 - s)
\end{array}
  \to J \to 0
$$
where $G = \bigoplus_{i=1}^s R(-e_i)$, $s \geq 0$, and $e_s > \ldots >
e_1 > e_0$.

Here, $s$ is the number of steps needed to reduce $J$ to the minimal
curve $J_0$ and $\beta_1, \beta_2, \beta_3 > 0$ are the Betti numbers
of $J_0$ (cf.\ Theorem 4.4).
\end{corollary}

\begin{proof}
The algorithm from \cite{MN7} that reduces the curves by using a facet of
maximal weight provides a sequence of reductions
$$
J=J_s \mapsto J_{s-1} \mapsto \cdots \mapsto J_1 \mapsto J_0
$$
where $J_0$ is a minimal curve.
Let $e_i$ be the maximal weight of a facet of the curve $J_i$ if $i >
0$ and let $e_0$ be the degree of the minimal generators of
$J_0$. Then Lemma 5.7 gives $e_s > \ldots > e_1 > e_0$ and the
resolution of $J$ is obtained by using Corollary 5.8 successively.
\end{proof}

A similar description can be given for the arithmetically
Cohen-Macaulay tetrahedral curves where we have to distinguish whether
the curve is componentwise linear or not. We leave the details to the
reader.

%%%%%%%%%%%%%%%%%%%%%%%%%%%%%%%%%%%%%%%%%%%%%%%%%%%%%%%%%%%%%%%%%%%%%%

\section{Tetrahedral Curves with Linear Resolutions} \label{linear resolutions}

Since the property of having a linear resolution is stronger  than
that of being componentwise linear, we now turn to the question of which
tetrahedral curves have linear resolution.  We begin with arithmetically
Cohen-Macaulay tetrahedral curves.

\begin{proposition} \label{acm lin res}
The following are the only arithmetically Cohen-Macaulay tetrahedral
curves with linear resolution (up to permutation of the variables):

\begin{itemize}
\item[(a)] $(r,0,0,0,0,0)$ for some $r \geq 1$;
\item[(b)] $(1,1,0,1,0,0)$ (this is the union of three non-coplanar lines
in ${\mathbb P}^3$ meeting at one point);
\item[(c)] $(1,1,1,1,1,1)$;
\item[(d)] $(2,1,0,1,0,1)$;
\item[(e)] $(2,1,1,1,1,2)$
\end{itemize}
\end{proposition}

\begin{proof}
It is easy to check that these curves do have linear resolution.  We have
to check that they are the only arithmetically Cohen-Macaulay tetrahedral
curves with this property (up to permutation of the variables).

Let $C$ be an arithmetically Cohen-Macaulay tetrahedral curve.  We
reduce using facets of maximal weight until one of the following happens:
either (i) we obtain a curve of type $(0,r,r,r,r,0)$ (up to permutation of
the variables), or (ii) we obtain a {\em plane} curve of degree 1, 2 or
3 (which is then one step away from the trivial curve via facets of
maximal weight).  In either of these cases, each step in the reduction,
passing from some
$J$ to a smaller curve with ideal  $I$, used a form $F$ that was not a
minimal generator of $I$, thanks to Theorem \ref{theorem B}.  It then
follows from Corollary \ref{first cor} and Theorem \ref{theorem B} (b)
that in case (i),
$I_C$ does not have a linear resolution.

So without loss of generality, we may reduce to a plane curve of degree
1, 2 or 3 via facets of maximal weight, each time using a form $F$ that is
not a minimal generator of the smaller curve.  If we arrive at a plane
curve of degree 2 or 3, then again by Corollary \ref{first cor}, $I_C$
does not have a linear resolution since a complete intersection of type
$(1,2)$ or type $(1,3)$ does not have a linear resolution.  So
$I_C$ reduces to a line via facets of maximal weight.

So we may work backwards, beginning with the curve $I = (1,0,0,0,0,0)$.
In order to form $J = L \cdot I + (F)$ and have the result have a linear
resolution, we need that $\deg F = 2$. Using (A) we obtain
$(2,0,0,0,0,0)$.  Using (C) we obtain $(1,1,0,1,0,0)$ (or
$(1,0,0,1,0,1)$ or $(1,1,0,0,0,1)$, which are equivalent).   (B) can only
repeat the result of (A), and (D) repeats the result of (C) (up to
permutation).

For the next step we have to pass from $(2,0,0,0,0,0)$ or $(1,1,0,1,0,0)$
to the next curve using $F$ of degree 3.  If we start with
$(2,0,0,0,0,0)$ and use (A) or (B), clearly $a_1$ becomes 3 so the
remaining entries must stay 0, and we can only obtain $(3,0,0,0,0,0)$.  If
we start with
$(2,0,0,0,0,0)$ and use (C) or (D) we obtain $(2,1,0,1,0,1)$ or 
$(2,0,1,0,1,1)$, which are
equivalent.  If we start with $(1,1,0,1,0,0)$ then the only permissible
basic double link that uses a form of degree 3 uses (D), and we obtain
$(1,1,1,1,1,1)$.

Passing to the next step, we need to use a form $F$ of degree 4.
Starting with $(3,0,0,0,$ $0,0)$, the only possibility is to use (A) or (B)
and pass to $(4,0,0,0,0,0,0)$.  If instead we start with $(2,1,0,1,0,1)$,
the only possibility is to use (D), from which we obtain
$(2,1,1,1,1,2)$.  If we start with $(1,1,1,1,1,1)$, none of the
operations produces a result with linear resolution since all forms $F$
will have degree 6.

For the next step, we need to use a form $F$ of degree 5.  From
$(4,0,0,0,0,0)$ it is clear that we can pass only to $(5,0,0,0,0,0)$.
 From $(2,1,1,1,1,2)$, none of the operations produces a result with
linear resolution.

It is clear that from $(r,0,0,0,0,0)$ we can obtain $(r+1,0,0,0,0,0)$.
Also, we know that if $I$ fails to have a linear resolution then so does
$J$, so once we lose this property we can never get it back.  Hence this
completes the proof.
\end{proof}

We now turn to non arithmetically Cohen-Macaulay tetrahedral curves.
Since basic double linkage preserves the even liaison class (\cite{LR},
\cite{migbook}), it is convenient to look within a fixed even liaison class.
Our first observation is that it can happen that there are fewer non-minimal
tetrahedral curves in the class than one might expect.

\begin{proposition} \label{no nonmin}
Let $I = (a_1,a_2,a_3,a_4,a_5,a_6)$ be a minimal tetrahedral curve.  Assume
that
\[
\begin{array}{rcl}
a_1 & > & \max \{ a_3 + a_5 +2,a_2+a_4 +2 \} \hbox{ and} \\
a_6 & > & \max \{ a_4+a_5 +2 ,a_2+a_3 +2 \}
\end{array}
\]
Then the even liaison class of $I$ contains no non-minimal tetrahedral
curves that reduce  to $I$.
\end{proposition}

\begin{proof}
It follows from Corollary \ref{s-min}.  If the even liaison class contained
a non-minimal tetrahedral curve that can be reduced via (A), (B), (C) and (D)
of Proposition \ref{MN7 basic result} to $I$, then in
the last step we pass from a tetrahedral curve $J$ to $I$, where the 6-tuple
corresponding to
$J$ has three of its entries equal to the corresponding ones of $I$, and
up to three others (and  exactly three others, if the entries are non-zero)
that are one more than the corresponding ones of $I$.  Without loss of
generality, suppose that  $J = (a_1+1, a_2+1, a_3+1,a_4,a_5,a_6)$.  But
the stated hypothesis then gives, via Corollary \ref{s-min}, that $J$ is
minimal.  Hence $J$ cannot have arisen from $I$ by basic double linkage.
\end{proof}

\begin{remark}
Since we do not yet have a good understanding of the Hartshorne-Rao module
of a tetrahedral curve, we do not know if there may be another 6-tuple that
is in the same even liaison class, also minimal, but which does allow
ascending tetrahedral curves.  Still, there are some cases where we know
that this does not happen.  For example, it was noted in \cite{MN7}, Remark
5.5, that the curve $(m,0,0,0,0,k)$, with $m,k \geq 2$, is the unique
minimal curve in its even liaison class, thanks to the main result of
\cite{LR}.   Since this curve satisfies the hypothesis of Proposition
\ref{no nonmin}, it is in fact the {\em only} tetrahedral curve in its even
liaison class.

We also remark that if a minimal tetrahedral curve admits one basic double
link of the type (A), (B), (C) or (D), then it allows infinitely many
(sequentially), and there are infinitely many tetrahedral curves in the
class.
\end{remark}

We have seen in Proposition \ref{acm lin res} that there are infinitely many
6-tuples representing arithmetically Cohen-Macaulay curves with linear
resolution, but if we identify those that are a multiple of a single line
then there are only finitely many.  We now show that the latter is true also
for non arithmetically Cohen-Macaulay curves.  We begin with the even
liaison class that we believe has the largest number of tetrahedral curves
with linear resolution.

\begin{proposition} \label{class of 2 lines}
Let $\mathcal L$ be the even liaison class of two skew lines.  Among
tetrahedral curves, this means $(1,0,0,0,0,1)$, $(0,1,0,0,1,0)$ or
$(0,0,1,1,0,0)$.  Then up to permutation of the variables, the following
are the only tetrahedral curves in $\mathcal L$ with linear resolution:

\begin{itemize}
\item[(a)] $(1,0,0,0,0,1)$

\item[(b)] $(2,1,0,0,0,1)$

\item[(c)] $(3,1,0,1,0, 1)$

\item[(d)] $(2,2,0,0,0,2)$

\item[(e)] $(2,1,1,1,0,1)$

\item[(f)] $(3,2,0,1,1,2)$

\item[(g)] $(3,2,1,1,2,3)$

\end{itemize}
\end{proposition}

\begin{proof}
Beginning with the ideal, $I$, of two skew lines, we must perform a basic
double link following the guidelines of Proposition \ref{MN7 basic result},
but using a form $F$ of degree 3 (since the generators of $I$ have degree
2).  So, for instance, from $(1,0,0,0,0,1)$ the only options are to use $F =
b^2c$ or $b^2d$ for type (A), $a^2c$ or $a^2d$ for type (B), $ad^2$ or
$bd^2$ for type (C), and
$ac^2$ or $bc^2$ for type (D), and we obtain the permutations of (b) having
either first or last entry equal to 2.  Continuing in this way (taking the
next basic double link using
$F$ of degree 4), one can exhaust all the possibilities.  We leave the
details to the reader.
\end{proof}

\begin{theorem} \label{finitely many curves}
Let $\mathcal L$ be the even liaison class of a non arithmetically
Cohen-Macaulay tetrahedral curve.  Then
$\mathcal L$ has only finitely many tetrahedral curves with linear
resolution.
\end{theorem}

\begin{proof}

Let $J$ be the ideal of a tetrahedral curve in ${\mathcal L}$ that has a
linear resolution.  We have seen that $J$ can be reduced to a minimal
tetrahedral curve by a sequence of reductions of the form (A), (B), (C) or
(D) as given in Proposition \ref{MN7 basic result}, and that we can do this
always using a facet of maximal weight (see also Definition
\ref{def of min}).  Let $C_0 = (a_1,a_2,a_3,a_4,a_5,a_6)$ be the minimal
curve so obtained.  We have seen in Theorem \ref{thm-res} that $I_{C_0}$ has
a linear resolution, and that the degree of its minimal generators is
$a_1+a_6$.

Our strategy will be to show that in any sequence of basic double links that
preserves the linearity of the resolutions, we cannot use any of (A), (B),
(C) or (D) more than once.

We first claim that in any such sequence of basic double links, all the
intermediate tetrahedral curves
$C_1, C_2,\dots$  between the ideals
$I_{C_0}$ and $J$ have linear resolution.  Indeed, suppose that $J$ reduces
to $I$, and that $I$ fails to have a linear resolution.  Assume that the
minimal generators of $J$ all have degree $d$. Since
$C_0$ is not arithmetically Cohen-Macaulay, we have seen (Theorem \ref{theorem
B}) that the polynomial $F$ used in the reduction is not a minimal generator
of $I$, but by construction it is a minimal generator of $J$ (it is the
only generator that does not have as a factor the linear form used in the
basic double link); hence it has degree $d$.  But then from the exact
sequence
\[
0 \rightarrow R(-d-1) \rightarrow I(-1) \oplus R(-d) \rightarrow J
\rightarrow 0
\]
it is clear that no splitting can occur in the mapping cone to restore
linearity to the non-linear resolution of $I$, no matter where the
non-linearity occurs in the resolution.

Consequently, if we work backwards, starting with $I_{C_0}$ and building up
to $J$ with basic double links, the first basic double link must use a
polynomial $F$ of degree $a_1+a_6+1$, the next a polynomial of degree
$a_1+a_6+2$, and so on.

If all entries $a_i > 0$, $1 \leq i \leq 6$, then the result is not hard to
see.  Indeed, suppose without loss of generality that the first basic double
link is of type (A) in Proposition \ref{MN7 basic result}.  Then $C_1 =
(a_1+1,a_2+1,a_3+1,a_4,a_5,a_6)$ and we have that the first three entries
give the facet of maximal weight.  Hence $a_1+a_2+a_3+3 = a_1+a_6+1$.  It is
clear that we cannot use type (A) again, since then the new curve $C_2$ will
have the sum of the first three entries be strictly greater than
$a_1+a_6+2$, while we would need equality.  But in fact, any other type that
we use increases one of $a_1,a_2,a_3$ by 1, so that we can {\em never}
return to use type (A).  But the same happens with the type used in the next
step -- it can be used at most once.  Continuing in this way, we see that at
most four basic double links can be used in order to preserve the linearity
of the resolution, so the result follows.

The only chance for the result to fail, then, is if some entries $a_i$ are
0, and remain 0 even after the basic double link (a possibility allowed in
Proposition \ref{MN7 basic result}).   Suppose without loss of generality
that $a_6 = \max \{ a_i \}$. We know from
Proposition \ref{MN7 basic result} and Corollary \ref{s-min} that also $a_1
>  0$, and that
\[
\begin{array}{rcl}
a_1 & > & \max \{ a_3 + a_5,a_2+a_4\} \hbox{ and} \\
a_6 & > & \max \{ a_4+a_5,a_2+a_3 \}.
\end{array}
\]
Suppose without
loss of generality, again, that $I_{C_1}$ is obtained via the basic double
link described in (A) of Proposition \ref{MN7 basic result}.  A priori,
$I_{C_1}$ could be any of the following tetrahedral curves:
\begin{itemize}
\item[(i)] $(a_1+1, 0 , 0 , a_4, a_5, a_6)$ (here $a_2=a_3=0$);
\item[(ii)] $(a_1+1, a_2+1, 0,a_4, a_5, a_6)$ (here $a_2 \geq 0, a_3=0$);
\item[(iii)] $(a_1+1, 0, a_3+1, a_4, a_5, a_6)$ (here $a_2 = 0, a_3 \geq 0$);
\item[(iv)] $(a_1+1,a_2+1,a_3+1,a_4, a_5, a_6)$;
\end{itemize}

\noindent Notice, though, that in case (i) the new curve $C_1$ is again
  a minimal curve (Corollary \ref{s-min}), so it is not in the same even
liaison class.  Hence (i) does not happen.  Also, in case (iv) it is easy to
see as above that we can never use (A) again, since the sum of the first
three entries is too big to use (A) to get $C_2$, and any of types (A)-(D)
increases at least one of $a_1,a_2$ or $a_3$ by 1 so $\deg F$ can never
``catch up" to get subsequent $C_i$.  By symmetry, there is no difference
between cases (ii) and (iii), so without loss of generality let us assume
that (ii) holds.

Since $\deg F = a_1+a_6+1 = (a_1+1)+(a_2+1)$ in the first basic double link,
we observe that
$a_6 = a_2 +1$, so
\[
I_{C_1} = (a_1+1, a_2+1, 0, a_4,a_5,a_2+1).
\]
Now, we have seen that the next basic double link (to pass from $C_1$ to
$C_2$) uses a polynomial
$F$ of degree $a_1+a_6+2$ in order to preserve linearity.  If we were to use
another basic double link of type (A), though, we would have to increase the
first entry and the second entry by 1.  Hence we would have
\[
\deg F = (a_1+2) + (a_2+2) + (0 \hbox{ or } 1) = a_1+a_6+2 = a_1+a_2 +3,
\]
which is impossible.  So we cannot use (A)
again, at least not now.

We will suppose that we use (B) for the second basic double link, and
carefully analyze the possibilities.  The other options for the second
basic double link are analyzed in a similar way.  If we use (B) for the
second basic double link, then we must  increase either the fourth entry or
the fifth entry (or both) by 1, since otherwise we have
$a_4+a_5=0
\geq a_6$, which is impossible.  Hence at least two entries (including $a_1$)
increase by 1, and as before, if all three entries increase by 1 then we can
never use (B) again.  So to preserve hope of using (B) again, we have two
cases: (i)
$a_4=0$ and remains 0 after the first application of (B), and (ii)
$a_5=0$ and remains 0 after the first application of (B).

In case (i), we have
\[
\begin{array}{rcl}
I_{C_1} & = & (a_1+1, a_2+1, 0, 0, a_5,a_2+1) \\
I_{C_2} & = & (a_1+2, a_2+1, 0, 0, a_5+1,a_2+1)
\end{array}
\]
But then we have
\[
(a_1+2) + 0 + (a_5+1) = a_1+a_6+2,
\]
so $a_6 = a_5+1$ and hence $a_2=a_5$.  Thus case (i) gives us
\[
\begin{array}{rcl}
I_{C_0} & = & (a_1,a_2, 0, 0, a_2,a_2+1)\\
I_{C_1} & = & (a_1+1, a_2+1, 0, 0, a_2,a_2+1) \\
I_{C_2} & = & (a_1+2, a_2+1, 0, 0, a_2+1,a_2+1)
\end{array}
\]
As before, we cannot use (B) again unless we use (C) at some point and
preserve $a_4=0$, since (A) and (D) both increase either $a_1$ or $a_5$ by
1.  And we cannot use (A) unless we use (D) first.

In case (ii) we have
\[
\begin{array}{rcl}
I_{C_1} & = & (a_1+1, a_2+1, 0, a_4, 0,a_2+1) \\
I_{C_2} & = & (a_1+2, a_2+1, 0, a_4+1,0,a_2+1)
\end{array}
\]
But then we have
\[
(a_1+2) + (a_4+1) +0 = a_1+a_6+2,
\]
so $a_6 = a_4+1$ and hence $a_2=a_4$.  Thus case (ii) gives
\[
\begin{array}{rcl}
I_{C_0} & = & (a_1,a_2, 0, a_2, 0,a_2+1)\\
I_{C_1} & = & (a_1+1, a_2+1, 0, a_2, 0,a_2+1) \\
I_{C_2} & = & (a_1+2, a_2+1, 0, a_2+1, 0,a_2+1)
\end{array}
\]
As before, we cannot use (B) again unless we use (D) at some point and
preserve $a_5=0$, since (A) and (C) both increase either $a_1$ or $a_4$ by
1.  And we cannot use (A) unless we use (D) first.

Now we consider the third basic double link (i.e. passing from $C_2$ to
$C_3$).  In case (i) above, we have two options: (i-a) to use (C) next, and
(i-b) to use (D) next.  In case (i-a), we obtain
\[
a_1+a_6+3 = (a_2+2) + (0 \hbox{ or } 1) + (a_2+2),
\]
which gives $a_1 = a_2 + (0 \hbox{ or } 1)$, and since $a_1 > a_2$ we have
$a_1 = a_2+1$.  But this means that our use of (C) increased $a_4$ from 0 to
1, and we can never use (B) again.  Similarly, since at this point the
second, fourth and sixth entries are $>0$, we can never use (C) again
either.  And at this stage we cannot use (A) again unless we use (D) and
preserve $a_3=0$. Hence
case (i-a) gives
\[
\begin{array}{rcl}
I_{C_0} & = & (a_2+1,a_2, 0, 0, a_2,a_2+1)\\
I_{C_1} & = & (a_2+2, a_2+1, 0, 0, a_2,a_2+1) \\
I_{C_2} & = & (a_2+3, a_2+1, 0, 0, a_2+1,a_2+1) \\
I_{C_3} & = & (a_2 +3, a_2+2, 0, 1, a_2+1, a_2+2)
\end{array}
\]
and the only possible fourth basic double link is (D).

In case (i-b) we have
\[
a_1+a_6+3 = (0 \hbox{ or } 1) + (a_2+2) + (a_2+2)
\]
so since $a_1 > a_2$, we have that (D) increases $a_3$ from 0 to 1, and $a_1
= a_2+1$.  In this case we can never use (A) or (D) again.
Hence case (i-b) gives
\[
\begin{array}{rcl}
I_{C_0} & = & (a_2+1,a_2, 0, 0, a_2,a_2+1)\\
I_{C_1} & = & (a_2+2, a_2+1, 0, 0, a_2,a_2+1) \\
I_{C_2} & = & (a_2+3, a_2+1, 0, 0, a_2+1,a_2+1) \\
I_{C_3} & = & (a_2 +3, a_2+1, 1, 0, a_2+2, a_2+2)
\end{array}
\]
Furthermore, since the fourth basic double link uses $\deg F = a_1+a_6+4$,
it is not hard to see that at this stage we cannot use (B) again.  Hence the
only possible fourth basic double link uses (C).

In both cases (i-a) and (i-b), it is not hard to see that the fourth basic
double link forces the last remaining 0 entry to become 1 (since we need
$\deg F = a_1+a_6+4 = 2a_2 + 6$), and hence we cannot use any of the four
types of basic double links and preserve the linearity of the resolution.
Case (ii), and the other cases, are proven similarly.
\end{proof}

%%%%%%%%%%%%%%%%%%%%%%%%%%%%%%%%%%%%%%%%%%%%%%%%%%%%%%%%%%%%%%%%%%%%%%

\section{The Generic Initial Ideal of a Tetrahedral Curve}

In this section, we will assume that the characteristic of $k$ is zero. With
this hypothesis, generic initial ideals are always strongly stable. We will
take generic initial ideals with respect to the reverse-lexicographic order;
using this order allows us to use some nice relationships from \cite{BS}
between an ideal and its gin.

Using Proposition \ref{ginres}, it is easy to describe the minimal 
generating set
of $\gin(J)$ when $J$ is the ideal of an arithmetically Cohen-Macaulay
tetrahedral curve.

\begin{proposition} \label{acmgin} Let $J$ be the ideal of an arithmetically
Cohen-Macaulay tetrahedral curve with lowest degree minimal generator in degree
$d_0$.

\begin{itemize}

\item[(a)] If $J$ is componentwise linear and has minimal generators in degrees
$d_0 \le \dots \le d_s$, then
$$
\gin(J) =
(a^{d_0},a^{d_0-1}b^{d_1-d_0+1},\dots,a^{d_0-p}b^{d_p-d_0+p},\dots,b^{d_s}).
$$
In particular, $s=d_0$, and $J$ has $d_0+1$ minimal generators.

\item[(b)] Suppose $J$ is not componentwise linear and that $J$ has $g$ minimal
generators in lowest degree $d_0$ and $h$ minimal generators in degree $d_0+1$.
Then $\gin(J)$ has minimal generating set with the monomials in $S=\{a^{d_0}$,
$a^{d_0-1}b$, $\dots, a^{d_0-(g-1)}b^{g-1}\}$, the first $h+g-1$ monomials of
degree $d_0+1$ not divisible by any element of $S$, and then elements of higher
degree. For each minimal generator of $J$ of degree higher than 
$d_0+1$, there is
a minimal generator $a^ib^j$ of $\gin(J)$ of the same degree with the powers on
$a$ decreasing down to zero.
\end{itemize}
\end{proposition}

\begin{proof} Since $J$ has codimension two, and $R/J$ is 
Cohen-Macaulay, $J$ and
$\gin(J)$ have projective dimension two. Because $\gin(J)$ is 
strongly stable, it
is generated by monomials of the form $a^ib^j$, including a pure power of both
$a$ and $b$; note that by the Eliahou-Kervaire resolution, any 
minimal generator
of $\gin(J)$ involving $c$ or $d$ would contradict the projective 
dimension being
two. Moreover, any stable ideal in two variables is a lexicographic ideal. It
follows immediately that
$$
\gin(J) = (a^{d_0}, a^{d_0-1}b^{q_1}, a^{d_0-2}b^{q_2}, \dots, b^{q_t}).
$$

Suppose first that $J$ is componentwise linear. Then $J$ and $\gin(J)$ have the
same graded Betti numbers and therefore minimal generators of the same degree.
The lowest degree generator of $\gin(J)$ is $a^{d_0}$. The second lowest has
degree $d_1$, and thus it is $a^{d_0-1}b^{d_1-(d_0-1)}$, and so on for the
others. All the generators have the form $a^{d_0-r}b^{d_r-(d_0-r)}$, 
where $0 \le
r \le s$. Note that the exponent on $a$ decreases by one as $r$ 
increases by one.
Since the minimal generator of $\gin(J)$ of highest degree is a pure power of
$b$, $s=d_0$, and $\gin(J)$ has $d_0+1$ minimal generators. Because $J$ and
$\gin(J)$ have the same graded Betti numbers, $J$ also has $d_0+1$ minimal
generators.

Part (b) follows from Proposition \ref{ginres}. The $(g-1)$ 
additional generators
in degree $d_0+1$ come from the fact that $(J_{d_0})$ does not have a linear
resolution, requiring us to add $g-1$ generators and syzygies of degree $d_0+1$
when we move to the gin.
\end{proof}

\begin{example} \label{ex:acmgin} Let $I$ be the ideal of the curve
$(1,2,2,2,1,2)$. Then $I$ has minimal resolution
$$
0 \rightarrow R(-7) \oplus R(-6)^2 \oplus R(-5) \rightarrow R(-6) \oplus
R(-5)^2 \oplus R(-4)^2 \rightarrow I \rightarrow 0.
$$
Note that $I$ is componentwise linear. Consequently,
$$
\gin(I)=(a^4,a^3b,a^2b^3,ab^4,b^6).
$$

Suppose now that $J$ is the ideal of the curve $(2,1,4,1,1,3)$. $J$ is not
componentwise linear, and it has minimal resolution
$$
0 \rightarrow R(-9) \oplus R(-8) \oplus R(-7)^2 \rightarrow R(-8) \oplus R(-7)
\oplus R(-6) \oplus R(-5)^2 \rightarrow J \rightarrow 0.
$$
Thus by Proposition \ref{acmgin}, $\gin(J)$ must have two minimal generators
of degree five and two minimal generators of degree six plus 
generators of higher
degree. Therefore
$$
\gin(J) = (a^5,a^4b,a^3b^3,a^2b^4,ab^6,b^8).
$$
\end{example}

We now turn to the generic initial ideal of a non arithmetically Cohen-Macaulay
tetrahedral curve.  We begin with a lemma that says that it is enough to
determine the generic initial ideal of the minimal curve in the even liaison
class.

\begin{lemma} \label{bdl-gin}
Let $J, I$ be the ideals of non-arithmetically Cohen-Macaulay tetrahedral
curves. Assume $J = L \cdot I + (F)$ is a basic double link of $I$ where $e :=
\deg F$ is the maximal weight of a facet of $J$. Then we have for the generic
initial ideals
$$
\gin (J) = a \gin (I) + b^e.
$$
\end{lemma}

\begin{proof}
By abuse of notation let us denote by $I$ and $J$ the ideals obtained from
$I$ and $J$ after a general change of coordinates. Then we have that $a \inn
(I) \subset \inn (J)$, hence $a \gin (I) \subset \gin (J)$. Since $\gin (J)$
is stable of codimension two, it must contain a power of $b$. We know by
\cite{BS} that the Castelnuovo-Mumford regularity of $J$ is $e$.
Therefore, we get
$$
a \gin(I) + (b^e) \subset \gin (J).
$$
But $a \gin (I) + (b^e)$ is a basic double link of $\gin (I)$. Since $I$ and
$J$ are componentwise linear, their graded Betti numbers agree with the
ones of their generic initial ideals. It follows that $a \gin (I) + (b^e)$
and $\gin (J)$ have identical graded Betti numbers, thus these ideals agree.
\end{proof}

While we are not yet able to determine the generic initial ideal of an
arbitrary minimal tetrahedral curve, we are able to do it for arithmetically
Buchsbaum tetrahedral curves.  It was shown in \cite{MN7} that up to a
permutation of variables, a minimal arithmetically Buchsbaum tetrahedral curve
is of the form $I_r = (r, 0, r-1, r-1, 0, r)$.  It is not hard to use
liaison addition (cf. \cite{schwartau}, \cite{GM4}) to show the recursive
relation
\begin{equation} \label{recursive}
I_{r+1} = (ac) \cdot I_{r} + (bd)^r \cdot I_1.
\end{equation}
(One shows the inclusion $\supseteq$ and then argues that the two ideals are
both saturated and define curves of the same degree.)

\begin{proposition}
The generic initial ideal of a minimal arithmetically Buchsbaum tetrahedral
curve $I_r = (r,0,r-1,r-1,0,r)$ is determined recursively by the following:
\begin{itemize}
\item[(a)] $\gin (I_1) = (a^2, ab, b^2, ac)$.

\item[(b)] $\gin (I_{r+1}) = (a^2) \cdot \gin (I_r) + (ab^{2r+1}, b^{2r+2},
a^{r+1}b^rc )$.

\end{itemize}

\end{proposition}

\begin{proof}
Part (a) is immediate, since $I_1$ has codimension two and is componentwise
linear, with four minimal generators all in degree 2, and
$\gin(I_1)$ is strongly stable.  For part (b), we have from (\ref{recursive})
that
\[
(a^2) \cdot \gin (I_r) \subseteq \gin(I_{r+1}).
\]
We also know that the number of minimal generators of $I_r$ is $3r+1$, all of
degree $2r$, and in fact that $I_r$ has a linear resolution:
\[
0 \rightarrow R(-2r-2)^r \rightarrow R(-2r-1)^{4r} \rightarrow R(-2r)^{3r+1}
\rightarrow I_r \rightarrow 0
\]
(cf.\ Theorem \ref{thm-res}).  Hence $\gin(I_r)$ has the same
resolution, since $I_r$ is componentwise linear.  From the above inclusion, we
have $3r+1$ minimal generators for $\gin(I_{r+1})$, and it is clear that also
$ab^{2r+1}$ and $b^{2r+2}$ are minimal generators, since $\gin(I_r)$ has
codimension two and is strongly stable.  We have only to prove that the last
minimal generator is $a^{r+1}b^rc$ (and not $a^{2r}c^2$, for instance).

Let $C := C_{r+1}$ be the tetrahedral curve with ideal $I_{r+1} = (r+1, 0,
r,r,0,r+1)$.  We know that $\deg C = 2(r+1)^2$, and that the
Hartshorne-Rao module $M(C)$ has dimension $r+1$ and is concentrated in
degree $2r$ (cf.\ \cite{BM1}).

Let $L$ be a general linear form defining a plane $H$ in ${\mathbb P}^3$
and let $t \in {\mathbb Z}$.  From the exact sequence
\[
0 \rightarrow (I_{C})_{t-1} \rightarrow (I_{C})_t \rightarrow
(I_{C \cap H})_t \rightarrow M(C)_{t-1} \rightarrow 0
\]
(where the last ``$\rightarrow 0$'' comes because $C = C_{r+1}$ is
arithmetically Buchsbaum), we see that
\[
\dim (I_{C\cap H})_t  =
\left \{
\begin{array}{ll}
0 & \hbox{ if $t \leq 2r$}; \\
r+1 & \hbox{ if $t = 2r+1$}.
\end{array}
\right.
\]
It follows that the $h$-vector of $C \cap H$ begins
\[
(1,2,3,\dots, 2r, 2r+1, r+1, \dots).
\]
But these entries already add up to $\deg C = 2(r+1)^2$, so this is the entire
$h$-vector.  It follows that in the quotient ring $S = R/(L) \cong k[a,b,c]$,
we have
\[
\gin(I_{C\cap H}) = (a^{2r}, a^{2r-1}b, \dots, a^{r+1}b^r, \dots)
\]
where the  entries up to $a^{r+1}b^r$ are all of degree $2r+1$ and the
remaining entries (not written) are of degree $2r+2$.

Since $I_C$ is saturated, without loss of generality we may reduce modulo $d$
and work in the ring $S = k[a,b,c]$.  Let $I = [I_{r+1}+(d)]/(d)$.  We will now
apply a result of
\cite{FG}, section 2.  They define (with our notation)
\[
I^j = \hbox{im} [(I:c^j) \rightarrow S \rightarrow S/(c)].
\]
We take $j=1$ and assume that we have a general change of coordinates, and
have taken the initial ideal.  Lemma 2.6 and Lemma 2.7 of \cite{FG} combine to
give that
\[
a^{i_0}b^{i_1} \in \gin(I^1) \iff a^{i_0} b^{i_1}c \in \gin(I).
\]
 From the information above about $C$ and $C\cap H$, it is clear that
$a^{r+1}b^r \in \gin (I^1)$; hence it follows that $a^{r+1}b^rc \in \gin(I)$,
and we have finished.
\end{proof}

%%%%%%%%%%%%%%%%%%%%%%%%%%%%%%%%%%%%%%%%%%%%%%%%%%%%%%%%%%%%%%%%%%%%%%


\begin{thebibliography}{999}

\bibitem {AAH} A.\ Aramova, J.\ Herzog, and T.\ Hibi, {\em Ideals with stable
Betti numbers}, \ Adv.\ Math. \ {\bf 152} (2000), no. 1, 72--77.

\bibitem{BS} D.\ Bayer and M.\ Stillman, {\em A criterion for detecting
{$m$}-regularity}, \ Invent.\ Math. \ {\bf 87} (1987), no. 1, 1--11.

\bibitem{BM1} G.\ Bolondi and J.\ Migliore, {\em Classification of Maximal Rank
Curves in the Liaison Class} $L_n$, Math.\ Ann.\ {\bf 277} (1987), 585--603.

\bibitem{ER} J.~A. Eagon and V.~Reiner, {\em Resolutions of
{S}tanley-{R}eisner rings and {A}lexander duality}, J.\ Pure Appl.\ Algebra
\textbf{130} (1998), no.~3, 265--275.

\bibitem{FG} G.\ Fl\o ystad and M.\ Green, {\em The information encoded in
initial ideals}, Trans.\ Amer.\ Math.\ Soc. {\bf 353} (2000), no.\ 4,
1427--1453.

\bibitem{Fran} C.~A. Francisco, {\em Resolutions of small sets of fat
points},  J.\ Pure Appl.\ Algebra (to appear).

\bibitem{GHP} V.~Gasharov, T.~Hibi, and I.~Peeva, {\em Resolutions of
{$\bold a$}-stable ideals}, J.\ Algebra \textbf{254} (2002), no.~2,
375--394.

\bibitem{GM4} A.V.\ Geramita and J.\ Migliore, {\em A Generalized Liaison
Addition}, J.\ Alg.\ 163 (1994), 139--164.

\bibitem{M2} D.~R. Grayson and M.~E. Stillman, \emph{Macaulay 2, a software
system for research in algebraic geometry}.
\newblock \verb|http://www.math.uiuc.edu/Macaulay2/|.

\bibitem{HH} J.~Herzog and T.~Hibi,  {\em Componentwise linear ideals},
Nagoya Math.\ J.\ \textbf{153} (1999), 141--153.

\bibitem{HRW} J.~Herzog, V.~Reiner, and V.~Welker, {\em Componentwise
linear ideals and Golod rings}, Michigan Math.\ J.\ \textbf{46} (1999), no.
2, 211--223.

\bibitem{LR} R.\ Lazarsfeld and P.\ Rao, {\em Linkage of General 
Curves of Large
Degree}, in ``Algebraic Geometry-- Open Problems (Ravello, 1982),'' Lecture
Notes in Mathematics, vol.\ 997, Springer--Verlag (1983), 267--289.

\bibitem{migbook} J.\ Migliore, ``Introduction to Liaison Theory and
Deficiency Modules,''  Birkh\"auser, Progress in Mathematics 165, 1998.

\bibitem{MN7} J.\ Migliore and U.\ Nagel, {\em Tetrahedral curves},
Int.\ Math.\ Res.\ Not.\ (to appear).

\bibitem{schwartau} P.\ Schwartau, {\em Liaison Addition and Monomial
Ideals},  Ph.D.\ thesis, Brandeis University (1982).

\end{thebibliography}
\end{document}